\documentclass{siamart220329}
\usepackage[all,cmtip]{xy}
\usepackage{graphicx}
\usepackage{float}   % For precise float placement
\usepackage{amssymb}
\usepackage{mathrsfs}
\usepackage{tikz-cd}
\usepackage{tikz}
\usepackage{caption}
\usepackage{subcaption}
\usepackage{enumitem}
\usetikzlibrary{fit, positioning, shapes.geometric}
\usepackage{setspace}
\usepackage{changepage}
\usepackage{euler}
\usepackage{leftidx}

\newsiamremark{example}{Example}
\newsiamremark{remark}{Remark}
\newsiamremark{notation}{Notation}
\crefname{example}{Example}{Examples}

\hypersetup{
  colorlinks = true,
  allcolors = siaminlinkcolor,
  urlcolor = siamexlinkcolor,
}

\headers{Short Title}{Authors}

\headers{Short Title}{Authors}

\title{Cellular sheaf Laplacians on the set of simplices of symmetric simplicial set induced by hypergraph}
\author{Seongjin Choi, Junyeong Park}
\date{November 2024}

\begin{document}

\maketitle

% REQUIRED
\begin{abstract}
We generalize cellular sheaf Laplacians on an ordered finite abstract simplicial complex to the set of simplices of a symmetric simplicial set. We construct a functor from the category of hypergraphs to the category of finite symmetric simplicial sets and define cellular sheaf Laplacians on the set of simplices of finite symmetric simplicial set induced by hypergraph. We provide formulas for cellular sheaf Laplacians and show that cellular sheaf Laplacian on an ordered finite abstract simplicial complex is exactly the ordered cellular sheaf Laplacian on the set of simplices induced by abstract simplicial complex.
\end{abstract}

% REQUIRED
\begin{keywords}
symmetric simplicial set, cellular sheaf, cellular sheaf cochain complex, cellular sheaf Laplacian, hypergraph
\end{keywords}

% REQUIRED
\begin{MSCcodes}
55U10, 55N05, 55N30
\end{MSCcodes}

\section{Introduction}
Graph Laplacian plays a central role in Graph Neural Networks which has wide real-world applications \cite{kipf2016semi, wu2019simplifying, buterez2024transfer, shlomi2020graph, dong2021adagnn}. For an abelian category $\mathscr{A}$, an $\mathscr{A}$-valued cellular sheaf $F$ on a graph $G$ \cite{curry2014sheaves} and a total order on $V(G)$, sheaf cochain complex of $G$ with $F$-coefficients $\delta : C^0(G, F) \to C^1(G, F)$ is defined \cite{hansen2019toward}. When $\mathcal{A}$ is the category of finite-dimensional $\mathbb{R}$-vector spaces, the adjoint of $\delta$, $\delta^*$, is well-defined and so is sheaf Laplacian $(\delta)^* \circ \delta$. Sheaf Laplacian is used to construct Sheaf Neural Networks \cite{hansen2020sheaf} on graphs for various purposes \cite{bodnar2022neural, barbero2022sheaf}.
Construction of sheaf Laplacian on the graph is generalized to an ordered finite abstract simplicial complex $L$ \cite{russold2022persistent}. $L$ is poset with respect to the set inclusion, so $L$ has a natural topology generated by Alexandrov base $\{U_{\sigma}:= \{\tau \in L \mid \sigma \subset \tau \} \}_{\sigma \in L}$. For simplex $\sigma$ of $L$, category $\mathscr{A}$ and $\mathscr{A}$-valued sheaf $\mathcal{F}$ on $L$, an associate cellular sheaf $F : L \to \mathscr{A}$ is a functor defined by $F(\sigma) := \mathcal{F}(U_{\sigma})$. Let $L_k$ be the set of $k$-simplices of $L$. Since any $\sigma \in L_k$ can be uniquely expressed as $\sigma = (v_0, \cdots, v_k)$ where $v_0 < \cdots < v_k$, for $l \in \{0, 1, \cdots, k \}$, face map $d_l : L_k \to L_{k-1}$ is defined by $d_l\left((v_0, \cdots, v_k)\right) := (v_0, \cdots,\widehat{v_l}, \cdots, v_k)$. When $\mathscr{A}$ is finitely complete abelian category, he defined cellular sheaf cochain complex of $L$ with $F$-coefficients, $\{C^k(L, F), \delta^k_F \}_{k \in \mathbb{Z}_{\geq 0}}$, given by
\begin{equation}\label{CSCCoforderedscpx}
C^k(L,F) := \bigoplus_{\sigma \in L_k} F(\sigma)    
\end{equation}
with the projection $\pi_{\sigma} : C^k(L,F) \to F(\sigma)$ and
\begin{equation}\label{coboundarymapoforderedscpx}
\delta^k_F := \bigoplus_{\tau \in L_{k+1}}\left(\sum_{l \in \{0, 1, \cdots, k \}} (-1)^l F(d_l \tau \subset \tau) \circ \pi_{d_l \tau} \right).    
\end{equation}
We call its cohomology as cellular sheaf cohomology. He showed that cellular sheaf cohomology is isomorphic to the sheaf cohomology of $L$ with $\mathcal{F}$-coefficients. As a corollary, when $\mathcal{A}$ is the category of finite-dimensional $\mathbb{R}$-vector spaces, degree $k$ cellular sheaf Laplacian $L^k_F$ is well-defined and its kernel is isomorphic to the degree $k$ sheaf cohomology of $L$ with $\mathcal{F}$-coefficients. Hence cellular sheaf Laplacian encodes information of both topology of $L$ and geometry of $F$. Cellular sheaf Laplacians for a trivial sheaf are used to signal processing on simplicial complexes \cite{barbarossa2020topological, yang2022simplicial, yang2023convolutional}.

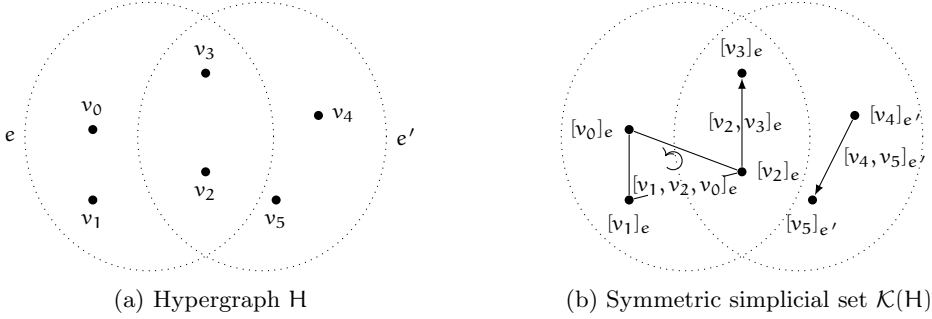
\begin{figure}[H]
    \centering
    \begin{subfigure}[b]{0.45\textwidth}
        \centering
        \begin{tikzpicture}
        [
            he/.style={draw, ellipse, dotted, inner sep=10pt},        % he = hyper edge
            node/.style={circle, fill=black, inner sep=0pt, minimum size=1.2mm, label distance=2mm} % node style with small black dot
        ]
        % hypergraph
        \node[node, label=above:{\footnotesize $v_0$}] (v0) at (-1.5,0.75) {};
        \node[node, label=below:{\footnotesize $v_1$}] (v1) at (-1.5,-0.1875) {};
        \node[node, label=below:{\footnotesize $v_2$}] (v2) at (0,0.1875) {};
        \node[node, label=above:{\footnotesize $v_3$}] (v3) at (0,1.5) {};
        \node[node, label=right:{\footnotesize $v_4$}] (v4) at (1.5,0.9375) {};
        \node[node, label=below:{\footnotesize $v_5$}] (g1) at (0.9375,-0.1875) {}; % to make e' large
        \node[draw=none, fill=none] (g2) at (0,0.84) {}; 

        \node[label=left:{\footnotesize $e$}] [he, fit=(v0) (v1) (v2) (v3)] {};
        \node[label=right:{\footnotesize $e'$}] [he, fit=(v2) (v3) (v4) (g1)] {};
        \end{tikzpicture}
        \caption{Hypergraph $H$}
    \end{subfigure}
    \hfill
    \begin{subfigure}[b]{0.45\textwidth}
        \centering
        \begin{tikzpicture}
        [
            he/.style={draw, ellipse, dotted, inner sep=10pt},        % he = hyper edge
            node/.style={circle, fill=black, inner sep=0pt, minimum size=1.2mm, label distance=2mm} % node style with small black dot
        ]
        % symm sset for hypergraph
        \node[node, label=left:{\footnotesize $[v_0]_e$}] (w0) at (4.5,0.75) {};
        \node[node, label=below:{\footnotesize $[v_1]_e$}] (w1) at (4.5,-0.1875) {};
        \node[node, label=right:{\footnotesize $[v_2]_e$}] (w2) at (6,0.1875) {};
        \node[node, label=above:{\footnotesize $[v_3]_e$}] (w3) at (6,1.5) {};
        \node[node, label=right:{\footnotesize $[v_4]_{e'}$}] (w4) at (7.5,0.9375) {};
        \node[node, label=below:{\footnotesize $[v_5]_{e'}$}] (a1) at (6.9375,-0.1875) {}; % to make e' large
        \node[draw=none, fill=none] (a2) at (4.875,-0.09375) {};
        \node[draw=none, fill=none] (a3) at (5.625,0.09375) {};

        \node[he, fit=(w0) (w1) (w2) (w3)] {};
        \node[he, fit=(w2) (w3) (w4) (a1)] {};
        % Draw directed edges between nodes with labels
        \draw[-latex] (w4) to node[midway, right] {\footnotesize $[v_4, v_5]_{e'}$} (a1);
        \draw[-latex] (w2) to node[midway, xshift=0.05cm] {\footnotesize $[v_2,v_3]_e$} (w3);
        \draw (w2) -- (w0) -- (w1) -- cycle;
        \draw (w1) -- (a2) -- cycle;
        \draw (a3) -- (w2) -- cycle;

        \draw[->] (5,0.25) arc[start angle=-135, end angle=150, radius=0.12];
        \node at (5.25,0) {\footnotesize $[v_1, v_2, v_0]_e$};
        \end{tikzpicture}
        \caption{Symmetric simplicial set $\mathcal{K}(H)$}
    \end{subfigure}
    \caption{(a) A hypergraph $H=(E(H), V(H), f_H)$ with $E(H) := \{e, e'\}$, $V(H) := \{ v_0, \cdots, v_5\}$, $f_H(e):= \{v_0, v_1, v_2, v_3\}$ and $f_H(e'):=\{v_2, v_3, v_4, v_5 \}$. (b) Description of symmetric simplicial set $\mathcal{K}(H)$ induced by $H$.}
    \label{figurehypergraph}
\end{figure}

We generalize cellular sheaf Laplacians on an ordered finite abstract simplicial complex to hypergraph. Since the hypergraph itself lacks mathematical structures, we associate finite symmetric simplicial set \cite{grandis2001finite} from hypergraph \cite{spivak2009higher} as in Figure \ref{figurehypergraph}. A symmetric simplicial set $X$, together with the set of simplices $\widehat{X}$ behave like an ordered finite abstract simplicial complex in the following sense : (1) each element of $\widehat{X}$ has dimension (2) $\widehat{X}$ has a natural preorder (3) there is a face map $d_l : X_n \to X_{n-1}$ for $l \in \{0, 1, \cdots, n\}$. These properties are sufficient to define $\mathscr{A}$-valued cellular sheaf $F$ on $\widehat{X}$ and cellular sheaf cochain complex of $\widehat{X}$ with $F$-coefficients without any choice of total order. We define the notion of cellular sheaf Laplacians of $\widehat{X}$ with $F$-coefficients when $\mathscr{A}$ is the category of finite dimensional $\mathbb{R}$-vector spaces. As in abstract simplicial complex, we show that the kernel of the degree $k$ cellular sheaf Laplacian is isomorphic to the degree $k$ sheaf cohomology when $X$ is closed, \v{C}ech.

\subsection{Organization}In section \ref{section2}, we review a cellular sheaf on a preordered set. We show that the category of cellular sheaves is equivalent to the category of sheaves. In section \ref{section3}, we review the symmetric simplicial set and its set of simplices. We define unordered, alternating, and ordered cellular sheaf cochain complexes on the set of simplices. We show that cellular sheaf cohomology is isomorphic to the sheaf cohomology when the symmetric simplicial set is closed, \v{C}ech. In section \ref{section4}, we construct a functor from the category of hypergraphs to the category of finite symmetric simplicial sets. We show that finite symmetric simplicial set induced by hypergraph is closed, \v{C}ech. When hypergraph is an ordered finite abstract simplicial complex, we show that the cellular sheaf cochain complex of an ordered finite abstract simplicial complex equals the ordered cellular sheaf cochain complex of the set of simplices induced by hypergraph. In section \ref{section5}, we compute explicit formulas of degree $k$ cellular sheaf Laplacians of a set of simplices induced by hypergraph.
\subsection{Conventions}
\begin{itemize}
    \item For categories $\mathscr{C}, \mathscr{D}$, we denote $[\mathscr{C}, \mathscr{D}]$ the category of functors from $\mathscr{C}$ to $\mathscr{D}$. We denote $\mathscr{C}(A,B):=\mathrm{Hom}_\mathscr{C}(A,B)$ for its hom-sets.
    \item We denote $\mathbf{Vect}_{\mathbb{R}}$ the category of finite dimensional inner product spaces over $\mathbb{R}$ whose morphisms are linear maps. $\mathbf{Vect}_{\mathbb{R}}$ is finitely bicomplete abelian category with enough injectives.
    \item We denote $\mathbf{Set}$ as the category of sets and $\mathbf{FinSet}$ as the category of finite sets. An endofunctor $\mathcal{P} : \mathbf{Set} \to \mathbf{Set}$ is defined by $\mathcal{P}(A) = \{ B \mid B \subseteq A, B \neq \emptyset \}$.
    \item For a topological space $X$ with base $\mathfrak{B}$, we denote $\mathrm{Sh}(X, \mathscr{A})$ the category of $\mathscr{A}$-valued sheaves on $X$ and $\mathrm{Sh}(\mathfrak{B}, \mathscr{A})$ the category of $\mathscr{A}$-valued $\mathfrak{B}$-sheaves.
    \item For a finite set $A$, we abbreviate $\{\underline{a}\} :=\{a_0, \cdots, a_n\} \subseteq A$ for a subset of $A$ and $(\underline{a}) := (a_0, \cdots, a_n) \in A^{n+1}$ for an element of $A^{n+1}$. We also abbreviate $(\underline{a})\backslash a_l := (a_0, \cdots, \hat{a}_l, \cdots, a_n)$ for $a_l \in \{\underline{a}\}$.
    \item For $n \in \mathbb{Z}_{\geq 0}$, we denote $\mathfrak{S}_n$ the set of bijections from $[n] := \{0, 1, \cdots, n\}$ to itself. For $g \in \mathfrak{S}_n$, we denote the sign of $g$ by $\mathrm{sgn}(g)$. For a set $A$, $n \in \mathbb{Z}_{\geq 0}$, $\mathfrak{S}_n$ acts on $A^{n+1}$ by $g \cdot (a_0, \cdots, a_n) := (a_{g(0)}, \cdots, a_{g(n)})$. We abbreviate $(\underline{a})_A \in \mathbf{FinSet}([n], A)$ a function such that $(\underline{a})_A(k) := a_k$ for $k \in [n]$.   
\end{itemize}

\section{Category of cellular sheaves on a preordered set}\label{section2}
Given a category $\mathscr{A}$ and a preordered set $P$, we define the category $\mathscr{A}$-valued cellular sheaves on a preordered set $P$ \cite{curry2018dualities}. When $\mathscr{A}$ is complete, we show that it is equivalent to the category of $\mathscr{A}$-valued sheaves on $P$.

\begin{definition}For a category $\mathscr{A}$, a preordered set $(P, \lesssim)$, the category $\mathrm{Cell}(P, \mathscr{A})$ is defined as follows.
\begin{itemize}
    \item The objects are all $F \in [P, \mathscr{A}]$ satisfying \[x \lesssim y, y \lesssim x \implies F(x)=F(y), F(y\lesssim x) = F(x \lesssim y) = \mathrm{Id}_{F(x)}. \] An object $F$ is called a $\mathscr{A}$-valued cellular sheaf on $P$.
    \item The morphisms from $F$ to $G$ are all natural transformations from $F$ to $G$.
\end{itemize}
\end{definition}
A preordered set has a natural topology, called the Alexandrov topology \cite{curry2018dualities}.
\begin{definition}Let $(P, \lesssim)$ be a preordered set.
\begin{itemize}
    \item For $p \in P$, \emph{the basic open set of $p$} is defined by $U_p := \{q \in P \ | \ q \gtrsim p \}$.
    \item $\mathfrak{P}:=\{U_p \}_{p \in P}$ forms a base for $P$. We call $\mathfrak{P}$ the \emph{Alexandrov base for $P$}.
    \item The topology generated by $\mathfrak{P}$ is called the \emph{Alexandrov topology on $P$}.
\end{itemize}
\end{definition}
The following lemma is a key property of the Alexandrov topology.
\begin{lemma}\label{Alexandrovtopology}Suppose $(P, \lesssim)$ is a preordered set equipped with Alexandrov topology. If an open set $U$ contains $p$, $U_p \subseteq U$.
\end{lemma}
\begin{proof} Suppose $U = \underset{q \in \Lambda}{\bigcup} U_q$ where $U_q \in \mathfrak{P}$ for each $q \in \Lambda \subset P$. Since $\{U_q\}_{q \in \Lambda}$ is an open cover of $U$, there exists $q_0 \in \Lambda$ satisfying $p \in U_{q_0}$. Hence $p \gtrsim q_0$ and $U_p \subseteq U_{q_0}\subseteq U$.
\end{proof}

Since we have Alexandrov topology on $P$, we can define the category of $\mathscr{A}$-valued sheaves on $P$, $\mathrm{Sh}(P, \mathscr{A})$.

\begin{proposition}\label{equivcellshv}Suppose $(P, \lesssim)$ is a preordered set and $\mathscr{A}$ is a complete category. Then $\mathrm{Cell}(P, \mathscr{A})$ and $\mathrm{Sh}(P, \mathscr{A})$ are equivalence of categories (See \cite[Proposition 3.3]{curry2018dualities} for poset $P$). 
\end{proposition}

\begin{proof} It suffices to show that $\mathrm{Cell}(P, \mathscr{A}) \cong \mathrm{Sh}(\mathfrak{P}, \mathscr{A})$ since $\iota : \mathrm{Sh}(\mathfrak{P}, \mathscr{A}) \to \mathrm{Sh}(P, \mathscr{A})$ defined by $\iota(\mathcal{F})(U) := \underset{U_p \subset U}{\varprojlim} \mathcal{F}(U_p)$ is an equivalence of categories which preserves basic open sets \cite[Tag 009H]{stacks-project}. Define a functor $\mathcal{S}' : \mathrm{Cell}(P, \mathscr{A}) \to \mathrm{PSh}(\mathfrak{P}, \mathscr{A})$ as 
\begin{itemize}
    \item $\mathcal{S}'(F)(U_p) := F(p)$ and for $p \lesssim q$, $res_{\mathcal{S}'(F)}(U_q \hookrightarrow U_p):=F(p \lesssim q)$.
    \item $\mathcal{S}'(\alpha)(U_p) := \alpha(p)$ for $\alpha \in \mathrm{Cell}(P, \mathscr{A})(F,G)$.
\end{itemize}
We show that $\mathcal{S}'(F)$ is a $\mathfrak{P}$-sheaf for any $F \in \mathrm{Cell}(P, \mathscr{A})$. Given $U_p \in \mathfrak{P}$, let $\{U_x\}_{x \in I}$ be a $\mathfrak{P}$-open cover. There is an element $x_0 \in I$ such that $U_{x_0}$ contains $p$, so $x_0 \lesssim p$. On the other hand, $U_{x_0} \subseteq U_{p}$ implies $p \lesssim x_0$. Hence $U_p = U_{x_0}$ and $\mathcal{S}'(F)(U_p)=\mathcal{S}'(F)(U_{x_0})=F(p)$. Given $x, y \in I$, let $\{U_{xyz} \}_{z \in I_{xy}}$ be a $\mathfrak{P}$-open cover of $U_x \cap U_y$. Suppose $\{s_x \in \mathcal{S}'(F)(U_x)\}_{x \in I}$ satisfies $s_x \vert_{U_{xyz}} = s_y \vert_{U_{xyz}}$. Then $s_{x_0} \in \mathcal{S}'(F)(U_p)$ is the unique section satisfying $s_{x_0} \vert_{U_x}=s_x$ for any $x \in I$. Hence $\mathcal{S}'(F)$ is a $\mathfrak{P}$-sheaf for any $F \in \mathrm{Cell}(P, \mathscr{A})$. \\
Let $\iota_P : P \to \mathfrak{P}^\mathrm{op}$ be a functor given by $\iota_P(q) := U_q$. Define a functor $\mathcal{T}' : \mathrm{Sh}(\mathfrak{P}, \mathscr{A}) \to \mathrm{Cell}(P, \mathscr{A})$ as
\begin{itemize}
    \item $\mathcal{T}'(\mathcal{F}) := \mathcal{F} \circ \iota_P$.
    \item $\mathcal{T}'(\eta)(p) := \eta(U_p)$ for $\eta \in \mathrm{Sh}(\mathfrak{P}, \mathscr{A})(\mathcal{F},\mathcal{G})$.
\end{itemize}
$(\mathcal{S}' \circ \mathcal{T}')(\mathcal{F})=\mathcal{F}$, $(\mathcal{T}' \circ \mathcal{S}')(F) = F$.\\
Hence $\mathrm{Cell}(P, \mathscr{A}) \cong \mathrm{Sh}(\mathfrak{P}, \mathscr{A})$ and $\mathcal{S}:= \iota \circ \mathcal{S}' : \mathrm{Cell}(P, \mathscr{A}) \to \mathrm{Sh}(P, \mathscr{A})$ is an equivalence of categories.
\end{proof}

\begin{remark}\label{remark} 
If $\{\mathcal{F}(U_p) \mid p \in P \}$ is finite set for any $\mathcal{F} \in \mathrm{Sh}(\mathfrak{P}, \mathscr{A})$, finitely completeness of $\mathscr{A}$ suffices to define $\iota$ by construction.
\end{remark}

\begin{definition}\label{sheafinducedbycellularsheaf}Let $(P, \lesssim)$ be a preordered set and $\mathscr{A}$ be a complete category. For $F \in \mathrm{Cell}(P, \mathscr{A})$, we say $\mathcal{S}(F)$ as \emph{the sheaf induced by $F$} where $\mathcal{S}$ is a functor in the proof of Proposition \ref{equivcellshv}. 
\end{definition}

\section{Cellular sheaf cochain complex on a set of simplices}\label{section3}
In this section, we develop a cellular sheaf theory on a symmetric simplicial set. More specifically, given a symmetric simplicial set $X$ and the set of simplices $\widehat{X}$, we define (1) a cellular sheaf $F$ on $\widehat{X}$ and sheaf $\mathcal{S}(F)$ on $\widehat{X}$ (2) cellular sheaf cochain complex of $\widehat{X}$ with $F$-coefficients (3) cellular sheaf Laplacians on $\widehat{X}$ for $\mathbf{Vect}_{\mathbb{R}}$-valued cellular sheaf $F$. We define conditions on $\widehat{X}$, called closed and \v{C}ech, to associate cellular sheaf cohomology with $F$-coefficients and sheaf cohomology with $\mathcal{S}(F)$-coefficients.

\subsection{Symmetric simplicial set and its set of simplices}In this subsection, we define symmetric simplicial sets \cite{grandis2001finite} and provide examples.

\begin{definition} A \emph{symmetric simplex category} $!\Delta$ is defined as follows.
\begin{itemize}
    \item The objects are $[n]$ for all $n\in\mathbb{Z}_{\geq 0}$.
    \item The morphisms from $[m]$ to $[n]$ are all functions from $[m]$ to $[n]$.
\end{itemize}
$[!\Delta^{\mathrm{op}}, \mathbf{Set}]$ is called the category of symmetric simplicial sets. 
\begin{itemize}
    \item An object $X$ of $[!\Delta^{\mathrm{op}}, \mathbf{Set}]$ is called a \emph{symmetric simplicial set}.
    \item A symmetric simplicial set $X$ is called \emph{finite} if $X \in [!\Delta^{\mathrm{op}}, \mathbf{FinSet}]$.
\end{itemize}
\end{definition}

\begin{definition}Let $X \in [!\Delta^{\mathrm{op}}, \mathbf{Set}]$, $n \in \mathbb{Z}_{\geq 0}$ and $l \in [n]$.
\begin{itemize}
    \item For $X \in [!\Delta^{\mathrm{op}}, \mathbf{Set}]$, we simply denote $X([n])$ as $X_n$. An element of $X_n$ is called a \emph{$n$-simplex} of $X$. A set $\widehat{X} := \underset{n \in \mathbb{Z}_{\geq 0}}{\coprod}X_n$ is called \emph{the set of simplicies of $X$}.
    \item A \emph{face map} $d_l : X_n \to X_{n-1}$ is defined by $d_l := X(d^l)$ where $d^l : [n-1] \to [n]$ is the function with $d^l(k) := k$ for $k < l$, $d^l(k) := k+1$ for $k \geq l$.
\end{itemize}
\end{definition}

Set of simplices of a symmetric simplicial set has a natural preorder.

\begin{proposition}\label{preorderonasetofsimpliciesofX} Suppose $X \in [!\Delta^{\mathrm{op}}, \mathbf{Set}]$ is a symmetric simplicial set. Define a relation $\lesssim$ on $\widehat{X}$ by $x \lesssim y$ if and only if there exists $\mu \in !\Delta([m], [n])$ satisfying $X(\mu)(y)=x$. Then $\lesssim$ becomes a preorder on $\widehat{X}$.
\end{proposition}
\begin{proof}$X(\mathrm{Id}_{[m]})(x)=x$, so $x \lesssim x$. Suppose $x \lesssim y$ and $y \lesssim z$. Then there are morphisms $\mu \in !\Delta([m], [n]), \nu \in !\Delta([n], [p])$ satisfying $X(\mu)(y)=x$ and $X(\nu)(z)=y$. Hence $X(\nu \circ \mu)(z)=\left(X(\mu) \circ X(\nu)\right)(z)=X(\mu)\left(X(\nu)(z)\right)=X(\mu)(y)=x$ and $x \lesssim z$.
\end{proof}

\begin{example} We introduce two classes of symmetric simplicial sets which will be used throughout the paper.
\begin{enumerate}[label=\textup{(\arabic*)}]
    \item Let $V$ be a set. A $V$-simplex $\Delta[V]$ is defined as follows.
        \begin{itemize}
            \item $\Delta[V]_n := \mathbf{Set}([n], V)$ for $n \in \mathbb{Z}_{\geq 0}$. $\Delta[V]_n$ is the set of all $(n+1)$-tuples in $V$.
            \item For $\mu : [m] \to [n]$ and $(\underline{v_i})_V \in \Delta[V]_n$, $(\Delta[V])(\mu) : \Delta[V]_n \to \Delta[V]_m$ is given by 
            \begin{align*}
                (\Delta[V])(\mu)\left((v_{i_0}, \cdots, v_{i_n})_V\right) &:= (v_{i_0}, \cdots, v_{i_n})_V \circ \mu \\
                &=(v_{i_{\mu(0)}}, \cdots, v_{i_{\mu(m)}})_V.
            \end{align*}
        \end{itemize}  
        Since $(\Delta[V])(\nu \circ \mu)= (\Delta[V])(\mu) \circ (\Delta[V])(\nu)$ for any $\mu : [m] \to [n]$, $\nu : [n] \to [p]$, $\Delta[V]$ is a symmetric simplicial set. When $V$ is a finite set, $\Delta[V]$ is a finite symmetric simplicial set.
    \item Let $X \in [!\Delta^{\mathrm{op}}, \mathbf{Set}]$. Since $\widehat{X}$ is a preordered set by Proposition \ref{preorderonasetofsimpliciesofX}, for $y \in X_0$, the basic open set $U_y$ is well-defined. A \v{C}ech nerve of $X$, denoted by $\check{C}(X)$, is defined as follows.
        \begin{itemize}
            \item $\check{C}(X)_n := \{ (\underline{y}):=(y_0, \cdots, y_n) \in (X_0)^{n+1} \mid U_{(\underline{y})}:=U_{y_0} \cap \cdots \cap U_{y_n} \neq \emptyset\}$
            \item For $\mu : [m] \to [n]$ and $(\underline{y_i}) \in \check{C}(X)_n$, $\check{C}(X)(\mu) : \check{C}(X)_n \to \check{C}(X)_m$ is given by
            \[\check{C}(X)(\mu)\left((y_{i_0}, \cdots, y_{i_n}) \right) := (y_{i_{\mu(0)}}, \cdots, y_{i_{\mu(m)}}).\]
        \end{itemize}
    For any $\mu : [m] \to [n]$, $\nu : [n] \to [p]$, $(\check{C}(X))(\nu \circ \mu)= (\check{C}(X))(\mu) \circ (\check{C}(X))(\nu)$. Hence $\check{C}(X)$ is a symmetric simplicial set. For $X, Y \in [!\Delta^{\mathrm{op}}, \mathbf{Set}]$ and $f \in [!\Delta^{\mathrm{op}}, \mathbf{Set}](X,Y)$, $\check{C}(f) : \check{C}(X) \to \check{C}(Y)$ is defined by \[\check{C}(f)_n\left((y_0, \cdots, y_n) \right):= (f_0(y_0), \cdots, f_0(y_n))\] for $(y_0 , \cdots, y_n) \in \check{C}(X)_n$. $\check{C}(f) \in [!\Delta^{\mathrm{op}}, \mathbf{Set}](\check{C}(X),\check{C}(Y))$ since for any $\mu : [m] \to [n]$ and $(y_0, \cdots, y_n) \in \check{C}(X)$, 
    \begin{align*}
    \check{C}(Y)(\mu) \circ \check{C}(f)_n\left((y_0, \cdots, y_n) \right) &= \check{C}(f)_m \circ \check{C}(X)(\mu)\left((y_0, \cdots, y_n) \right)\\
    &= \left(f_0(y_{\mu(0)}), \cdots, f_0(y_{\mu(n)}) \right).
    \end{align*}
    Suppose $X, Y, Z \in [!\Delta^{\mathrm{op}}, \mathbf{Set}]$, $f \in [!\Delta^{\mathrm{op}}, \mathbf{Set}](X,Y)$ and $g \in [!\Delta^{\mathrm{op}}, \mathbf{Set}](Y,Z)$. Then $\check{C}(g \circ f) = \check{C}(g) \circ \check{C}(f)$, so $\check{C} : [!\Delta^{\mathrm{op}}, \mathbf{Set}] \to [!\Delta^{\mathrm{op}}, \mathbf{Set}]$ is a functor. 
       
\end{enumerate}
\end{example}

We define two properties of a symmetric simplicial set, called closed and \v{C}ech.

\begin{definition}\label{Cechdefn}Let $X \in [!\Delta^{\mathrm{op}}, \mathbf{Set}]$. For $n \in \mathbb{Z}_{\geq 0}$, we have a map $\psi_n : X_n \to \check{C}(X)_n$ given by $\psi_n(y) := \Bigl(X((0)_{[n]})(y), \cdots, X((n)_{[n]})(y) \Bigr)$. We say $X$ is
\begin{enumerate}[label=\textup{(\arabic*)}]
    \item \emph{closed} if $\{\emptyset \} \cup \{U_y \}_{y \in \widehat{X}}$ is closed under finite intersections.
    \item \emph{\v{C}ech} if $\psi =\{\psi_n\} : X \to \check{C}(X)$ is an isomorphism in $[!\Delta^{\mathrm{op}}, \mathbf{Set}]$.
\end{enumerate}
\end{definition}

\begin{remark}\label{cechfunctoriality}Suppose $X, Y \in [!\Delta^{\mathrm{op}}, \mathbf{Set}]$ are \v{C}ech and $f \in [!\Delta^{\mathrm{op}}, \mathbf{Set}](X,Y)$. Then $\check{C}(f) \circ \psi = \psi \circ f$ since for any $n \in \mathbb{Z}_{\geq 0}$, $y \in X_n$, 
\begin{align*}
\check{C}(f)_n \circ \psi_n(y) &=\left(f_0(X((0)_{[n]})(y)), \cdots, f_0(X((n)_{[n]})(y)) \right) \\
&= \left(Y((0)_{[n]})(f_n(y)), \cdots, Y((n)_{[n]})(f_n(y)) \right)(\because \ f \in [!\Delta^{\mathrm{op}}, \mathbf{Set}](X,Y))\\
&=\psi_n \circ f_n(y). 
\end{align*}
In particular, if $f$ is an isomorphism, so is $\check{C}(f)$.    
\end{remark}

\subsection{\v{C}ech cochain complexes} In this subsection, we recall the notion of \v{C}ech cochain complexes \cite{gallier2022homology}.

\begin{definition}Let $Y$ be a topological space with an open cover $\mathcal{W}=\{W_{\beta}\}_{\beta \in B }$ of $Y$. For an abelian category $\mathscr{A}$, let $\mathcal{F}$ be an $\mathscr{A}$-valued presheaf on $Y$. Let $k \in \mathbb{Z}_{\geq 0}$.
\begin{itemize}
    \item For $(\underline{\beta})\in B^{k+1}$, we define $W_{(\underline{\beta})}:=W_{\beta_0}\cap \cdots \cap W_{\beta_k}$.
    \item We define \[\check{C}^k(\mathcal{W}, \mathcal{F}):=\underset{(\underline{\beta})\in B^{k+1}}{\prod}\mathcal{F}(W_{(\underline{\beta})})\] \emph{the set of all unordered $k$-cochains of $\mathcal{W}$ with $\mathcal{F}$-coefficients}. An element of $\check{C}^k(\mathcal{W}, \mathcal{F})$ is called an \emph{unordered $k$-cochain of $\mathcal{W}$ with $\mathcal{F}$-coefficients}. We abbreviate a \emph{$k$-cochain} when $\mathcal{W}$ and $\mathcal{F}$ are clear. We denote the $(\underline{\beta})$-projection by $\pi_{(\underline{\beta})} : \check{C}^k(\mathcal{W}, \mathcal{F}) \to \mathcal{F}(W_{(\underline{\beta})})$.
    \item \emph{A coboundary map} $\delta^k_{\mathcal{W}, \mathcal{F}} : \check{C}^k(\mathcal{W}, \mathcal{F}) \to \check{C}^{k+1}(\mathcal{W}, \mathcal{F})$ is defined as \[\delta^k_{\mathcal{W}, \mathcal{F}} := \underset{(\underline{\beta})}{\prod}\left(\sum_{l\in [k+1]}(-1)^l\mathrm{res}_{\mathcal{F}}(W_{(\underline{\beta})} \hookrightarrow{} W_{(\underline{\beta}\backslash \beta_l)}) \circ \pi_{(\underline{\beta})\backslash \beta_l}\right).\] Computation shows that $\delta^{k+1}_{\mathcal{W}, \mathcal{F}} \circ \delta^k_{\mathcal{W}, \mathcal{F}} = 0$. $\{\check{C}^k(\mathcal{W}, \mathcal{F}), \delta^k_{\mathcal{W}, \mathcal{F}}\}_{k \in \mathbb{Z}_{\geq 0}}$ is called an \emph{unordered \v{C}ech cochain complex of $\mathcal{W}$ with $\mathcal{F}$-coefficients}.
    \item A $k$-cochain $s$ is called \emph{alternating} if for any $g \in \mathfrak{S}_{k}$, \[\pi_{g \cdot (\beta_0, \cdots, \beta_k) } \circ s = \pi_{(\beta_{g(0)}, \cdots, \beta_{g(k)})} \circ s = \mathrm{sgn}(g)\cdot \pi_{(\beta_0, \cdots, \beta_k)} \circ s.\] We denote \emph{the set of all alternating $k$-cochains} by $\check{C}^k_{\mathrm{alt}}(\mathcal{W}, \mathcal{F})$. $\delta^k_{\mathcal{W}, \mathcal{F}}$ induces a map $\delta^k_{\mathcal{W}, \mathcal{F}} : \check{C}^k_{\mathrm{alt}}(\mathcal{W}, \mathcal{F}) \to \check{C}^{k+1}_{\mathrm{alt}}(\mathcal{W}, \mathcal{F})$. $\{\check{C}^k_{\mathrm{alt}}(\mathcal{W}, \mathcal{F}), \delta^k_{\mathcal{W}, \mathcal{F}}\}_{k \in \mathbb{Z}_{\geq 0}}$ is called \emph{an alternating \v{C}ech cochain complex of $\mathcal{W}$ with $\mathcal{F}$-coefficients}.
    \item Suppose $(B, <)$ is a totally ordered set. We define \[\check{C}^k_{\mathrm{ord}}(\mathcal{W}, \mathcal{F}):=\underset{(\beta_0 < \cdots < \beta_k)\in B^{k+1}}{\prod}\mathcal{F}(W_{(\beta_0, \cdots, \beta_k)})\] \emph{the set of all ordered $k$-cochains of $\mathcal{W}$ with $\mathcal{F}$-coefficients}. An element of $\check{C}^k_{\mathrm{ord}}(\mathcal{W}, \mathcal{F})$ is called an \emph{ordered $k$-cochain of $\mathcal{W}$ with $\mathcal{F}$-coefficients}. $\delta^k_{\mathcal{W}, \mathcal{F}}$ induces a map $\delta^k_{\mathcal{W}, \mathcal{F}} : \check{C}^k_{\mathrm{ord}}(\mathcal{W}, \mathcal{F}) \to \check{C}^{k+1}_{\mathrm{ord}}(\mathcal{W}, \mathcal{F})$. $\{\check{C}^k_{\mathrm{ord}}(\mathcal{W}, \mathcal{F}), \delta^k_{\mathcal{W}, \mathcal{F}}\}_{k \in \mathbb{Z}_{\geq 0}}$ is called an \emph{ordered \v{C}ech cochain complex of $\mathcal{W}$ with $\mathcal{F}$-coefficients}.
    \item Suppose $(B,<)$ is a totally ordered set. For $q \in \mathbb{Z}_{\geq 0}$, we denote
    \begin{itemize}
        \item $\check{H}^q((C^k(\mathcal{W}, \mathcal{F}), \delta^k_{\mathcal{W}, \mathcal{F}}))$ the \emph{$q$th unordered \v{C}ech cohomology of $\mathcal{W}$ with $\mathcal{F}$-coefficients}
        \item $\check{H}^q((\check{C}^k_{\mathrm{alt}}(\mathcal{W}, \mathcal{F}), \delta^k_{\mathcal{W}, \mathcal{F}}))$ the \emph{$q$th alternating \v{C}ech cohomology of $\mathcal{W}$ with $\mathcal{F}$-coefficients}
        \item $\check{H}^q((\check{C}^k_{\mathrm{ord}}(\mathcal{W}, \mathcal{F}), \delta^k_{\mathcal{W}, \mathcal{F}}))$ the \emph{$q$th ordered \v{C}ech cohomology of $\mathcal{W}$ with $\mathcal{F}$-coefficients}.
    \end{itemize}
    Three cohomologies are all isomorphic \cite[Tag 01FG]{stacks-project}.
    \item $\mathrm{Cover}(Y)$ is the category of open covers of $Y$ such that
    \begin{itemize}
        \item the objects are open covers of $Y$ and
        \item the morphisms from $\mathcal{U}$ to $\mathcal{V}$ is $\{*\}$ when $\mathcal{V}$ is a refinement of $\mathcal{U}$ and otherwise $\emptyset$. 
    \end{itemize}
    \item 
    The $q$th unordered \v{C}ech cohomology of $Y$ with $\mathcal{F}$-coefficients as \[\check{H}^q(Y, \mathcal{F}) := \underset{\underset{\mathcal{W} \in \mathrm{Cover}(Y)}{\longrightarrow}}{\lim} \check{H}^q(\mathcal{W}, \mathcal{F}).\]
    The $q$th alternating \v{C}ech cohomology of $Y$ with $\mathcal{F}$-coefficients is defined as \[\check{H}^q_{\mathrm{alt}}(Y, \mathcal{F}) := \underset{\underset{\mathcal{W} \in \mathrm{Cover}(Y)}{\longrightarrow}}{\lim} \check{H}^q_{\mathrm{alt}}(\mathcal{W}, \mathcal{F}).\]
    The $q$th ordered \v{C}ech cohomology of $Y$ with $\mathcal{F}$-coefficients is defined as \[\check{H}^q_{\mathrm{ord}}(Y, \mathcal{F}) := \underset{\underset{\mathcal{W}=\{W_{\beta} \}_{\beta \in (B, <)} \in \mathrm{Cover}(Y)}{\longrightarrow}}{\lim} \check{H}^q_{\mathrm{ord}}(\mathcal{W}, \mathcal{F}).\]
\end{itemize}
\end{definition}

\subsection{Cellular sheaf cochain complex}Let $X$ be a finite symmetric simplicial set, $\mathscr{A}$ be an abelian category and $F \in \mathrm{Cell}(\widehat{X}, \mathscr{A})$ be a $\mathscr{A}$-valued cellular sheaf on $\widehat{X}$. In this subsection, we define three versions of cellular sheaf cochain complexes of $\widehat{X}$ with $F$-coefficients.

\begin{definition}Let $X \in [\Delta^{\mathrm{op}}, \mathbf{FinSet}]$, $\mathscr{A}$ be an abelian category, $F \in \mathrm{Cell}(\widehat{X}, \mathscr{A})$ and $k \in \mathbb{Z}_{\geq 0}$.
\begin{itemize}
 \item We define $C^k(\widehat{X}, F) := \underset{y \in X_k}{\bigoplus}F(y)$ \emph{the set of all unordered cellular sheaf $k$-cochains of $\widehat{X}$ with $F$-coefficients}. An element of $C^k(\widehat{X}, F)$ is called an \emph{unordered $k$-cochain of $\widehat{X}$ with $F$-coefficients}. We abbreviate a cellular sheaf $k$-cochain when $X$ and $F$ are clear. We denote $\pi_{y, F} : C^k(\widehat{X}, F) \to F(y)$ the $y$-projection.
 \item \emph{A coboundary map} $\delta^k_F : C^k(\widehat{X}, F) \to C^{k+1}(\widehat{X}, F)$ is defined as
\begin{equation}\label{cellularcoboundarymap}
\delta^k_F := \underset{{z \in X_{k+1}}}{\bigoplus} \left(\sum_{l \in [k+1]}(-1)^l \cdot F(d_l(z)\lesssim z) \circ \pi_{d_l(z)}\right).
\end{equation}
Computation shows that $\delta^{k+1}_F \circ \delta^k_F =0$. $\{C^k(\widehat{X},F), \delta^k_F\}_{k \in \mathbb{Z}_{\geq 0}}$ is called an \emph{unordered cellular sheaf cochain complex of $\widehat{X}$ with $F$-coefficients}.
\item A cellular sheaf $k$-cochain $s$ is called \emph{alternating} if for any $g \in \mathfrak{S}_{k}$,
\[\pi_{X(g)(y)} \circ s = \mathrm{sgn}(g)\cdot \pi_{y} \circ s.\]
We denote \emph{the set of all alternating cellular sheaf $k$-cochains} by $C^k_{\mathrm{alt}}(\widehat{X},F)$. $\delta^k_F$ induces a map $\delta^k_F : C^k_{\mathrm{alt}}(\widehat{X},F) \to C^{k+1}_{\mathrm{alt}}(\widehat{X},F)$. $\{C^k_{\mathrm{alt}}(\widehat{X},F), \delta^k_F\}_{k \in \mathbb{Z}_{\geq 0}}$ is called an \emph{alternating cellular sheaf cochain complex of $\widehat{X}$ with coefficients in $F$}. 
\item Suppose $X$ is \v{C}ech and $(X_0,<)$ is a totally ordered set. We define 
\begin{equation}
C^k_{\mathrm{ord}}(\widehat{X}, F) := \underset{(\underline{y})=(y_0 < \cdots < y_k) \in \check{C}(X)_k}{\bigoplus}F(\psi^{-1}_k(\underline{y}))
\end{equation}
\emph{the set of all ordered cellular sheaf $k$-cochains of $\widehat{X}$ with $F$-coefficients}. An element of $C^k_{\mathrm{ord}}(\widehat{X}, F)$ is called an \emph{ordered $k$-cochain of $\widehat{X}$ with $F$-coefficients}. $\delta^k_F$ induces a map $\delta^k_F : C^k_{\mathrm{ord}}(\widehat{X},F) \to C^{k+1}_{\mathrm{ord}}(\widehat{X},F)$. $\{C^k_{\mathrm{ord}}(\widehat{X},F), \delta^k_F\}_{k \in \mathbb{Z}_{\geq 0}}$ is called an \emph{ordered cellular sheaf cochain complex of $\widehat{X}$ with $F$-coefficients}.
\item $H^q((C^k(\widehat{X},F), \delta^k_F))$ is called the \emph{$q$th unordered cellular sheaf cohomology of $\widehat{X}$ with $F$-coefficients}. $H^q((C^k_{\mathrm{alt}}(\widehat{X},F), \delta^k_F))$ is called the \emph{$q$th alternating cellular sheaf cohomology of $\widehat{X}$ with $F$-coefficients}. When $X$ is \v{C}ech and $(X_0, <)$ is totally ordered, $H^q((C^k_{\mathrm{ord}}(\widehat{X},F), \delta^k_F))$ is called the \emph{$q$th ordered cellular sheaf cohomology of $\widehat{X}$ with $F$-coefficients}.
\end{itemize}
\end{definition}

\begin{lemma}\label{cechCSCCisomorphism}Suppose $X, Y$ are finite symmetric simplicial sets and $f:= \{f_n \}_{n \in \mathbb{Z}_{\geq 0}}  : X \to Y$ is an isomorphism in $[!\Delta^{\mathrm{op}}, \mathbf{FinSet}]$. Suppose $F$ is an $\mathscr{A}$-valued cellular sheaf on $\widehat{X}$ for an abelian category $\mathscr{A}$. Define $f_* F \in \mathrm{Cell}(\widehat{Y}, \mathscr{A})$ by $(f_*F)(y):= F(f^{-1}_n(y))$ for $y \in Y_n$. 
\begin{enumerate}[label=\textup{(\arabic*)}]
    \item A bijection $\widehat{f} : \widehat{X} \to \widehat{Y}$ defined by $\widehat{f}(y) := f_k(y)$ for $y \in X_k$ is homeomorphism and the following diagram
    \begin{align*}
    \xymatrix{
    \mathrm{Cell}(\widehat{X}, \mathscr{A}) \ar[d]_-{\mathcal{S}} \ar[r]^-{f_*} & \mathrm{Cell}(\widehat{Y}, \mathscr{A}) \ar[d]^-{\mathcal{S}} \\
    \mathrm{Sh}(\widehat{X}, \mathscr{A}) \ar[r]_-{\widehat{f}_*} & \mathrm{Sh}(\widehat{Y}, \mathscr{A})
    }
\end{align*}
is commutative. Moreover, $\mathcal{S} = (\widehat{f}^{-1})_* \circ \mathcal{S} \circ f_*$.
    \item $(C^k(\widehat{X}, F), \delta^k_F) = (C^k(\widehat{Y}, f_*F), \delta^k_{f_*F})$ and $(C^k_{\mathrm{alt}}(\widehat{X}, F), \delta^k_F) = (C^k_{\mathrm{alt}}(\widehat{Y}, f_*F), \delta^k_{f_*F})$ for any $k \in \mathbb{Z}_{\geq 0}$.
    \item Suppose $X,Y$ are \v{C}ech, $(X_0, <) , (Y_0, \prec)$ are totally ordered sets and $f_0$ preserves the total order. Then $(C^k_{\mathrm{ord}}(\widehat{X}, F), \delta^k_F) = (C^k_{\mathrm{ord}}(\widehat{Y}, f_*F), \delta^k_{f_*F})$ for any $k \in \mathbb{Z}_{\geq 0}$.
\end{enumerate}
\end{lemma}

\begin{proof}(1) For any $z \in Y_k$, $\widehat{f}^{-1}(U_z)=U_{f_k^{-1}(z)}$ is open in $\widehat{X}$. Hence $\widehat{f}$ is continuous. Same argument shows that $\widehat{f}^{-1}$ is continuous. For any open set $U \in \widehat{Y}$ and $F \in \mathrm{Cell}(\widehat{X}, \mathscr{A})$,
\begin{align*}
\left(\mathcal{S}(f_*F) \right)(U) &= \underset{q \in U}{\varprojlim}(f_*F)(q) \\
&=\underset{q \in U}{\varprojlim}F\left(\widehat{f}^{-1}(q) \right) \\
&=\underset{p \in \widehat{f}^{-1}(U)}{\varprojlim}F(p) \\
&=\mathcal{S}(F)\left(\widehat{f}^{-1}(U) \right)\\
&=\widehat{f}_*\left(\mathcal{S}(F) \right)(U).
\end{align*}
Hence $\mathcal{S} \circ f_* = \widehat{f}_* \circ \mathcal{S}$. Functoriality of the direct image functor proves $\mathcal{S} = (\widehat{f}^{-1})_* \circ \mathcal{S} \circ f_*$.

(2) Direct computations show that
\begin{align*}
C^k(\widehat{Y}, f_*F) &= \underset{z \in Y_k}{\bigoplus}\left(f_*F\right)(z)=\underset{z \in Y_k}{\bigoplus}F(f^{-1}_k(z)) \\
&= \underset{y \in X_k}{\bigoplus}F(y)(\because \textrm{$f_k$ is bijective})\\
&=C^k(\widehat{X}, F).
\end{align*}
Since $f$ is an isomorphism, for any $k \in \mathbb{Z}_{\geq 0}, g \in \mathfrak{S}_k$, $Y(g) \circ f_k = f_k \circ X(g)$ and $f_k$ is bijective. Hence
\begin{align*}
&s \in C^k_{\mathrm{alt}}(\widehat{X},F) \\
&\iff \pi_{X(g)(y), F} \circ s = \mathrm{sgn}(g)\cdot \pi_{y, F} \circ s  \textrm{ for any $g \in \mathfrak{S}_k, y \in X_k$} \\
&\iff \pi_{f_k^{-1}(Y(g)(f_k(y))), F} \circ s = \mathrm{sgn}(g) \cdot \pi_{f_k^{-1}(f_k(y)), F} \circ s\\
&\iff \pi_{f_k^{-1}(Y(g)(z)), F} \circ s = \mathrm{sgn}(g) \cdot \pi_{f_k^{-1}(z), F} \circ s \textrm{ for any $g \in \mathfrak{S}_k$, $z \in Y_k$}\\
&\iff \pi_{(Y(g)(z)), f_*F} \circ s = \mathrm{sgn}(g) \cdot \pi_{z, f_*F} \circ s \textrm{ for any $g \in \mathfrak{S}_k, z \in Y_k$} \\
&\iff s \in C^k_{\mathrm{alt}}(\widehat{Y},f_*F).
\end{align*}

Since $f$ is an isomorphism, for any $k \in \mathbb{Z}_{\geq 0}$ and $l \in [k]$, $f^{-1}_k \circ d_l = d_l \circ f^{-1}_{k+1}$. Hence

\begin{align*}
\delta^k_{f_* F} &= \underset{z \in Y_{k+1}}{\bigoplus}\left(\sum_{l \in [k+1]}(-1)^l \cdot (f_*F)(d_l(z) \lesssim z) \circ \pi_{d_l(z)}\right)\\   
&= \underset{z \in Y_{k+1}}{\bigoplus}\left(\sum_{l \in [k+1]}(-1)^l \cdot F\left(f^{-1}_k(d_l(z)) \lesssim f^{-1}_{k+1}(z)\right) \circ \pi_{f^{-1}_k(d_l(z))}\right) \\
&= \underset{z \in Y_{k+1}}{\bigoplus}\left(\sum_{l \in [k+1]}(-1)^l \cdot F\left(d_l(f^{-1}_{k+1}(z)) \lesssim f^{-1}_{k+1}(z)\right) \circ \pi_{d_l(f^{-1}_{k+1}(z))}\right) \\
&=\underset{y \in X_{k+1}}{\bigoplus}\left(\sum_{l \in [k+1]}(-1)^l \cdot F\left(d_l(y) \lesssim y\right) \circ \pi_{d_l(y)}\right)\\
&=\delta^k_F.
\end{align*}

(3) Since $f$ is an isomorphism, $\check{C}(f_k) \circ \psi_k = \psi_k \circ f_k$, $f_k^{-1} \circ \psi_k^{-1} = \psi^{-1}_k \circ \left(\check{C}(X)(f_k)\right)^{-1}$ and $\check{C}(X)(f_k)$ is an isomorphism by Remark \ref{cechfunctoriality}. Hence
\begin{align*}
C^k_{\mathrm{ord}}(\widehat{Y}, f_*F) &= \underset{(\underline{z})=(z_0 \prec \cdots \prec z_k) \in \check{C}(Y)_k}{\bigoplus}(f_*F)(\psi^{-1}_k((\underline{z})))\\
&= \underset{(\underline{z})=(z_0 \prec \cdots \prec z_k) \in \check{C}(Y)_k}{\bigoplus}F(f_k^{-1} \circ \psi^{-1}_k((\underline{z})))\\
&= \underset{(\underline{z})=(z_0 \prec \cdots \prec z_k) \in \check{C}(Y)_k}{\bigoplus}F\left(\psi^{-1}_k \circ \left(\check{C}(X)(f_k)\right)^{-1}((\underline{z}))\right)\\
&= \underset{(\underline{y})=(y_0 < \cdots < y_k) \in \check{C}(X)_k}{\bigoplus}F\left(\psi^{-1}_k((\underline{y}))\right)\\
&=C^k_{\mathrm{ord}}(\widehat{X}, F).
\end{align*}
\end{proof}

\subsection{Isomorphisms of cohomologies}
Let $Y$ be a topological space and $\mathscr{A}$ be an abelian category with enough injectives. For $\mathcal{F} \in \mathrm{Sh}(Y, \mathscr{A})$ and $q \in \mathbb{Z}_{\geq 0}$, we denote $H^q_{\mathrm{sh}}(Y, \mathcal{F})$ the $q$th sheaf cohomology of $Y$ with $\mathcal{F}$-coefficients. Suppose $X$ is \v{C}ech finite symmetric simplicial set and $\mathscr{A}$ is a complete abelian category with enough injectives. Then $F \in \mathrm{Cell}(\widehat{X}, \mathscr{A})$ induces $\psi_* F \in \mathrm{Cell}(\widehat{\check{C}(X)},\mathscr{A})$ and $\mathcal{S}(\psi_*F) \in \mathrm{Sh}(\widehat{\check{C}(X)}, \mathscr{A})$. We have three cohomologies : (1) $H^q_{\mathrm{sh}}(\widehat{\check{C}(X)}, \mathcal{S}(\psi_*F)) \cong H^q_{\mathrm{sh}}(\widehat{X}, (\widehat{\psi}^{-1})_*(\mathcal{S}(\psi_*F)))\cong H^q_{\mathrm{sh}}(\widehat{X}, \mathcal{S}(F))$ by Lemma \ref{cechCSCCisomorphism} (2) $\check{H}^q(\widehat{\check{C}(X)}, \mathcal{S}(\psi_*F))$ (3) $H^q(\widehat{X}, F) \cong H^q(\widehat{\check{C}(X)}, \psi_*F)$ by Lemma \ref{cechCSCCisomorphism}. We show three cohomologies are isomorphic when $X$ is closed, \v{C}ech.
\begin{theorem}\label{Cechsheaf}
Suppose $X \in [!\Delta^{\mathrm{op}}, \mathbf{FinSet}]$ is a finite symmetric simplicial set and $\mathscr{A}$ is a complete abelian category with enough injectives. Suppose $F \in \mathrm{Cell}(\widehat{X}, \mathscr{A})$ is an $\mathscr{A}$-valued cellular sheaf on the set of simplicies of $X$ and $\mathcal{F} \in \mathrm{Sh}(\widehat{X}, \mathscr{A})$ is an $\mathscr{A}$-valued sheaf on the set of simplicies of $X$.
\begin{enumerate}[label=\textup{(\arabic*)}]
\item If $X$ is closed, $H^q_{\mathrm{sh}}(\widehat{X}, \mathcal{F}) \cong \check{H}^q(\widehat{X}, \mathcal{F})$ for $q \in \mathbb{Z}_{\geq 0}$.
\item If $X$ is \v{C}ech with isomorphism $\psi : X \to \check{C}(X)$ in Definition \ref{Cechdefn} and $X_0$ is totally ordered, \[\check{H}^q(\widehat{\check{C}(X)}, \mathcal{S}(\psi_*F)) \cong H^q(\widehat{X}, F) \cong H^q_{\mathrm{alt}}(\widehat{X}, F) \cong  H^q_{\mathrm{ord}}(\widehat{X}, F)\] for $q \in \mathbb{Z}_{\geq 0}$. \\
\item If $X$ is closed and \v{C}ech, \[H^q_{\mathrm{sh}}(\widehat{X}, \mathcal{S}(F)) \cong H^q(\widehat{X}, F) \cong H^q_{\mathrm{alt}}(\widehat{X}, F) \cong  H^q_{\mathrm{ord}}(\widehat{X}, F)\] for $q \in \mathbb{Z}_{\geq 0}$.
\end{enumerate}
\end{theorem}

\begin{proof}
(1) It suffices to show that the Alexandrov base $\mathfrak{X}=\{\emptyset \cup  U_y\}_{y \in \widehat{X}}$ for $\widehat{X}$ satisfies (a) $\mathfrak{X}$ is closed under finite intersections (b) $\check{H}^k_{\mathrm{alt}}(U_{y_0}\cap \cdots \cap U_{y_k}, \mathcal{F})=0$ for any $y_0, \cdots, y_k \in \widehat{X}$ and $k>0$ by the Cartan's theorem \cite[Theorem 13.19.]{gallier2022homology}. (a) is satisfied since $X$ is closed. Closedness of $X$ implies that when $U_{y_0}\cap \cdots \cap U_{y_k}$ is nonempty, it should be $U_y$ for some $y \in \widehat{X}$. Hence we will show that $\check{H}^k_{\mathrm{alt}}(U_y, \mathcal{F})=0$ to prove (b). \\
Suppose $\mathcal{W}=\{W_{\beta} \}_{\beta \in B}$ is an open cover of $U_y$. There exists $\beta_0 \in B$ such that $y \in W_{\beta_0}$, $U_y \subset W_{\beta_0}$ by Lemma \ref{Alexandrovtopology}. Hence $r : \{\ast\} \to B$ with $r(\ast) := \beta_0$ implies $\{ U_y\}_{\{ \ast \} }$ refines $\mathcal{W}$. Since $\mathcal{W}$ was arbitrary, $\{ U_y\}_{\{ \ast \} }$ is a terminal object in $\mathrm{Cover}(U_y)$. Therefore, $\check{H}^k_{\mathrm{alt}}(U_y, \mathcal{F}) \cong \check{H}^k_{\mathrm{alt}}(\{U_y\}_{\{ \ast \} }, \mathcal{F})=0$ for $k>0$ due to the fact that $\check{C}^k_{\mathrm{alt}}(\{ U_y\}_{\{ \ast \} }, \mathcal{F})=0$ for $k >0$. \\
(2) Consider a collection of open sets $\mathcal{W}_{\mathrm{term}, X} = \{U_v \}_{v \in X_0}$ in $\widehat{X}$. For any $y \in X_n$, $X\left((0)_{[n]} \right)(y)  \lesssim y$ for $X\left((0)_{[n]} \right)(y) \in X_0$ and $y \in U_{X\left((0)_{[n]} \right)(y)}$. Since $y$ was arbitrary, $\mathcal{W}_{\mathrm{term}, X}$ covers $\widehat{X}$ and $\mathcal{W}_{\mathrm{term}, X} \in \mathrm{Cover}(\widehat{X})$.\\
Suppose $\mathcal{W}=\{W_{\beta} \}_{\beta \in B} \in \mathrm{Cover}(\widehat{X})$. Given $v \in X_0$, there exists $\beta_v \in B$ satisfying $v \in W_{\beta_v}$. Define a map $r : X_0 \to B$ defined by $r(v) := \beta_v$. Since $\beta_v \in B$, $U_v \subset W_{\beta_v}$ by Lemma \ref{Alexandrovtopology} and $\mathcal{W}_{\mathrm{term},X}$ refines $\mathcal{W}$. Since $\mathcal{W}$ was arbitrary, $\mathcal{W}_{\mathrm{term},X}$ is a terminal object in $\mathrm{Cover}(\widehat{X})$. Hence the \v{C}ech cochain complex of $\widehat{\check{C}(X)}$ is isomorphic to the \v{C}ech cochain complex of $\mathcal{W}_{\mathrm{term},\check{C}(X)}$ and

\begin{align*}
\check{C}^k(\mathcal{W}_{\mathrm{term},\check{C}(X)}, \mathcal{S}(\psi_*F)) &= \bigoplus_{(\underline{y_i}) \in \check{C}(X)_k}\mathcal{S}(\psi_*F)(U_{(\underline{y_i})})(\because \textrm{$\mathscr{A}$ is abelian category}) \\
&= \bigoplus_{(\underline{y_i}) \in \check{C}(X)_k}(\psi_*F)\left((\underline{y_i})\right)(\because \textrm{$X$ is \v{C}ech and Proposition \ref{equivcellshv}}) \\
&=C^k(\widehat{\check{C}(X)}, \psi_*F)\\
&= C^k(\widehat{X}, F)(\because \textrm{Lemma \ref{cechCSCCisomorphism}}).
\end{align*}
Lemma \ref{cechCSCCisomorphism} also implies \[\check{C}^k_{\mathrm{alt}}(\mathcal{W}_{\mathrm{term},\check{C}(X)}, \mathcal{S}(\psi_*F)) = C^k_{\mathrm{alt}}(\widehat{X},F)\] and \[\check{C}^k_{\mathrm{ord}}(\mathcal{W}_{\mathrm{term},\check{C}(X)}, \mathcal{S}(\psi_*F)) = C^k_{\mathrm{ord}}(\widehat{X},F).\]
Coboundary map $\delta^k_{\mathcal{W}_{\mathrm{term},\check{C}(X)}, \mathcal{S}(\psi_*F)}$ is given by
\begin{align*}
&\underset{(\underline{y_j}) \in \check{C}(X)_{k+1}}{\bigoplus}\left(\sum_{l \in [k+1]} (-1)^l \mathrm{res}_{\mathcal{S}(\psi_*F)}(U_{(\underline{y_j})} \hookrightarrow{} U_{(\underline{y_j})\backslash y_l}) \circ \pi_{(\underline{y_j})\backslash y_l}\right) \\   
&= \underset{z \in X_{k+1}}{\bigoplus}\left(\sum_{l \in [k+1]}(-1)^l\cdot F\left(d_l(z) \lesssim z\right) \circ \pi_{d_l(z)}\right) \\
&= \delta^k_F.
\end{align*}
Hence 
\[\check{H}^q(\widehat{\check{C}(X)}, \mathcal{S}(\psi_*F)) \cong H^q(\widehat{X}, F) \cong H^q_{\mathrm{alt}}(\widehat{X}, F) \cong  H^q_{\mathrm{ord}}(\widehat{X}, F).\]
(3)
\begin{align*}
H^q_{\mathrm{sh}}(\widehat{X}, \mathcal{S}(F)) &\cong H^q_{\mathrm{sh}}(\widehat{X}, (\widehat{\psi}^{-1})_*(\mathcal{S}(\psi_*F)))(\because \textrm{Lemma \ref{cechCSCCisomorphism}})\\
&\cong H^q_{\mathrm{sh}}(\widehat{\check{C}(X)}, \mathcal{S}(\psi_*F))(\because \textrm{$\widehat{\psi}$ is homeomorphism})\\
&\cong \check{H}^q(\widehat{\check{C}(X)}, \mathcal{S}(\psi_*F))(\because \textrm{Theorem \ref{Cechsheaf}.(1)})\\
&\cong H^q(\widehat{\check{C}(X)}, \psi_*F)(\because \textrm{Theorem \ref{Cechsheaf}.(2)})\\
&\cong H^q(\widehat{X}, F)((\because \textrm{Lemma \ref{cechCSCCisomorphism}})\\
&\cong H^q_{\mathrm{alt}}(\widehat{X}, F) \cong  H^q_{\mathrm{ord}}(\widehat{X}, F)(\because \textrm{Theorem \ref{Cechsheaf}.(2)}).
\end{align*}

\end{proof}

\subsection{Cellular sheaf Laplacians}Let $X$ be a finite simplicial set and $F$ is a $\mathbf{Vect}_{\mathbb{R}}$-valued cellular sheaf $F$ on $\widehat{X}$. In this subsection, we define degree $k$ cellular sheaf Laplacians on $\widehat{X}$ for $k \in \mathbb{Z}_{\geq 0}$.

\begin{definition}\label{cellularsheaflaplacian} Suppose $X \in [\Delta^{\mathrm{op}}, \mathbf{FinSet}]$, $F \in \mathrm{Cell}(\widehat{X}, \mathbf{Vect}_{\mathbb{R}})$ and $k \in \mathbb{Z}_{\geq 0}$.
\begin{itemize}
    \item Let $\langle, \rangle_{C^k}$ be the induced inner product on $C^k(\widehat{X}, F)$. $(\delta^k_F)^* : C^{k+1}(\widehat{X},F) \to C^k(\widehat{X},F)$ is called \emph{the adjoint of $\delta^k_F$} if 
    \begin{equation}\label{adjointmap}
    \langle \delta^k_F(s) , s' \rangle_{C^{k+1}} = \langle s , (\delta^k_F)^*(s') \rangle_{C^{k}}     
    \end{equation}
    for any $s \in C^{k}(\widehat{X},F), s' \in C^{k+1}(\widehat{X},F)$.
    \item The \emph{degree $k$ up-Laplacian of $F$ on $\widehat{X}$} is a linear map $L^k_{F,+} : C^k(\widehat{X},F) \to C^k(\widehat{X},F)$ defined by $L^k_{F,+}:=(\delta^k_F)^*\delta^k_F$. The \emph{degree $k$ down-Laplacian of $F$ on $\widehat{X}$} is a linear map $L^k_{F,-} : C^k(\widehat{X},F) \to C^k(\widehat{X},F)$ defined by $L^k_{F,-}:=\delta^{k-1}_F(\delta^{k-1}_F)^*$. Set $L^0_{F,-} := 0$. The \emph{degree $k$ cellular sheaf Laplacian of $F$ on $\widehat{X}$} is a linear map $L^k_F : C^k(\widehat{X},F) \to C^k(\widehat{X},F)$ defined by $L^k_F := L^k_{F,+} + L^k_{F,-}$.
    \item The restriction of $L^k_{F, +}, L^k_{F, -}, L^k_{F}$ to $C^k_{\mathrm{alt}}(\widehat{X}, F)$ are called \emph{degree $k$ alternating up-Laplacian, down-Laplacian and Laplacian of $F$ on $\widehat{X}$}. We denote by $L^k_{\mathrm{alt}, F, +}, L^k_{\mathrm{alt}, F, -}$ and $L^k_{\mathrm{alt}, F}$.
    \item Suppose $\check{X}$ is \v{C}ech and $(X_0, <)$ is totally ordered set. The restriction of $L^k_{F, +}, L^k_{F, -}, L^k_{F}$ to $C^k_{\mathrm{ord}}(\widehat{X}, F)$ are called \emph{degree $k$ ordered up-Laplacian, down-Laplacian and Laplacian of $F$ on $\widehat{X}$}. We denote by $L^k_{\mathrm{ord}, F, +}, L^k_{\mathrm{ord}, F, -}$ and $L^k_{\mathrm{ord}, F}$.
\end{itemize}
\end{definition}

\begin{lemma}\label{adjointkerimagelemma}Suppose $f : V \to W$ is a linear map between inner product spaces and $f^* : W \to V$ is the adjoint of $f$.
    \begin{enumerate}[label=\textup{(\arabic*)}]
        \item $\mathrm{Ker} \ f \cap \mathrm{Im} \ f^*= \{0_V\}$.
        \item $\mathrm{Ker} \ f^* \cap \mathrm{Im} \ f= \{0_W\}$.
    \end{enumerate}
\end{lemma}

\begin{proof}
(1) For $v \in \mathrm{Ker} \ f \cap \mathrm{Im} \ f^*$, $v = f^*(w)$ for some $w \in W$ and $\langle v, v \rangle = \langle v, f^*(w) \rangle = \langle f(v), w \rangle = \langle 0_V, w \rangle = 0_V$. Hence $v = 0$. \\
(2) For $w \in \mathrm{Ker} \ f^* \cap \mathrm{Im} \ f$, $w=f(v)$ for some $v \in V$ and $\langle w , w \rangle = \langle f(v), w \rangle = \langle v, f^*(w) \rangle = \langle v , 0_W \rangle = 0$. Hence $w = 0_W$.
\end{proof}

We show that $L^k_F$ contains topological information about $\widehat{X}$ when $\widehat{X}$ is closed and \v{C}ech.
\begin{theorem}Suppose $X \in [\Delta^{\mathrm{op}}, \mathbf{FinSet}]$ is a fnite symmetric simplicial set which is closed and \v{C}ech. Suppose $F \in \mathrm{Cell}(\widehat{X}, \mathbf{Vect}_{\mathbb{R}})$ is a cellular sheaf on $\widehat{X}$ such that $\widehat{X}$ $\{F(z) \mid z \in \widehat{X} \}$ is finite set.
Then $\mathrm{Ker} \ L^k_F \cong H^k_{\mathrm{sh}}(\widehat{X}, \mathcal{S}(F))$.
\end{theorem}

\begin{proof}It suffices to show that $\mathrm{Ker} \ L^k_F \cong H^k(\widehat{X}, F)$ by Remark \ref{remark} and Theorem \ref{Cechsheaf}.(3). $\mathrm{Ker} \ L^k_F = \mathrm{Ker} \ \delta^k_F \cap \mathrm{Ker} \ (\delta^{k-1}_F)^*$ since 
\begin{align*}
\langle L^k_F(s), s \rangle &= \langle (\delta^k_F)^* \delta^k_F(s) , s \rangle + \langle \delta^{k-1}_F(\delta^{k-1}_F)^*(s), s\rangle\\
&= \langle \delta^k_F(s), \delta^k_F(s)\rangle + \langle (\delta^{k-1}_F)^*(s), (\delta^{k-1}_F)^*(s)\rangle\\
&= \Vert \delta^k_F(s) \Vert^2 + \Vert (\delta^{k-1}_F)^*(s) \Vert^2.
\end{align*}
Hence $L^k_F(s)=0$ if and only if $\delta^k_F(s)=0, (\delta^{k-1}_F)^*(s)=0$.

Choose a basis of $\mathrm{Ker} \ L^k_F$ and extend to $C^k(\widehat{X}, F)$. Then $C^k(\widehat{X}, F) = \mathrm{Ker} \ L^k_F \oplus \mathrm{Im} \ L^k_F$. Define a linear map $f : \mathrm{Ker} \ L^k_F \to H^k(\widehat{X}, F)$ by $f(s) := [s]$. $f$ is well-defined since $\mathrm{Ker} \ L^k_F = \mathrm{Ker} \ \delta^k_F \cap \mathrm{Ker} \ (\delta^{k-1}_F)^* \subset \mathrm{Ker} \ \delta^k_F$. To show $f$ is surjective, suppose $[s] \in H^k(\widehat{X}, F)$ for some $s \in C^k(\widehat{X},F)$ with $\delta^k_F(s)=0$. Then $s = s' + L^k_F(s'')$ for some $s' \in \mathrm{Ker} \ L^k_F$, $s'' \in C^k(\widehat{X},F)$. Since $s, s' \in \mathrm{Ker} \  \delta^k_F$, $\delta^k_F(L^k_F(s''))=\delta^k_F(\delta^k_F)^*\delta^k_F(s'')=0$ and $(\delta^k_F)^*\delta^k_F(s'') \in \mathrm{Ker} \ \delta^k_F \cap \mathrm{Im} \ (\delta^k_F)^* = \{0\}$ by Lemma \ref{adjointkerimagelemma}. Hence $(\delta^k_F)^*\delta^k_F(s'')=0$ and $s=s' + (\delta^{k-1}_F) \circ (\delta^{k-1}_F)^*(s'')$. This implies $f(s')=[s]$.
To show $f$ is injective, suppose $f(s) =[s]=0$ for some $s \in \mathrm{Ker} \ L^k_F$. Then $s= \delta^{k-1}_F(s')$ for some $s' \in C^k(\widehat{X},F)$. Since $s \in \mathrm{Ker} \ L^k_F = \mathrm{Ker} \ \delta^k_F \cap \mathrm{Ker} \ (\delta^{k-1}_F)^*$, $s \in \mathrm{Ker} \ (\delta^{k-1}_F)^* \cap \mathrm{Im}\ \delta^{k-1}_F = \{0\}$ by Lemma \ref{adjointkerimagelemma}. Hence $s = 0$ and $f$ is injective.
\end{proof}
\section{Finite symmetric simplicial set induced by hypergraph}\label{section4}
In this section, we define a functor from the category of hypergraphs to the category of finite symmetric simplicial sets. Since an ordered finite abstract simplicial complex $L$ is also hypergraph, we have a finite symmetric simplicial set induced by $L$. We show that cellular sheaf cochain complex of $L$ is the ordered cellular sheaf cochain complex of set of simplices induced by $L$.

\subsection{Category of hypergraphs}
\begin{definition}Let $T : \mathbf{FinSet} \to \mathbf{FinSet}$ be an endofunctor.
\begin{itemize}
    \item A \emph{$T$-graph} $X= \left(E(X), V(X), f_X : E(X) \to T(V(X))\right)$ is an object of the comma category $(\mathrm{Id} \downarrow T)$ \cite{jakel2015coalgebraic}. $E(X)$ is called the \emph{edge set of $X$}, $V(X)$ is called the \emph{vertex set of $X$} and $f_X$ is called the \emph{structure map of $X$}.
    \item A $\mathcal{P}$-graph $H=(E(H),V(H), f_H)$ is called an \emph{hypergraph} if $f_H(e) \notin V(H)$ for any $e \in E(H)$. We denote $\mathscr{H}$ \emph{the category of hypergraphs}. 
    \item For a hypergraph $H=(E(H),V(H), f_H)$, the \emph{extended structure map} $\tilde{f}_H : E(H) \coprod V(H) \to \mathcal{P}(V(H))$ is defined by $\tilde{f}_H((0,e)) :=f_H(e), \tilde{f}_H((1,v)):=v$ for $e \in E(H), v \in V(H)$.
\end{itemize}  
\end{definition}

\begin{example}\label{examplehypergraph} 
\begin{enumerate}[label=\textup{(\arabic*)}]
    \item Suppose $L$ is a finite abstract simplicial complex. Then $L$ induces a natural hypergraph $(E(L), V(L), f_L) = (L \backslash L_0, L_0, \mathrm{Id}_{L \backslash L_0})$. We abuse the notation $L$ for indicating hypergraph $(L \backslash L_0, L_0, \mathrm{Id}_{L \backslash L_0})$.
    \item A set of triples $H=(E(H), V(H), f_H)$ given by 
        \begin{itemize}
            \item $E(H) := \{e, e'\}$
            \item $V(H) := \{ v_0, \cdots, v_5\}$
            \item $f_H(e):= \{v_0, v_1, v_2, v_3\}$ and $f_H(e'):=\{v_2, v_3, v_4, v_5 \}$
        \end{itemize}
        is hypergraph. Figure \ref{figurehypergraph}.(a) describes the geometric description of $H$.
\end{enumerate}
\end{example}

\subsection{Construction of a functor}In this subsection, we construct a functor $\mathcal{K} : \mathscr{H} \to [!\Delta^{\mathrm{op}}, \mathbf{FinSet}]$ inspired by \cite{spivak2009higher}. Geometrically, $(\mathcal{K}(H))_n$ is the disjoint unions of all $n$-simplices of $\Delta[\tilde{f}_H(x)]$ with identifying all same $(n+1)$-tuples as in Figure \ref{figurehypergraph}.(b).

\begin{theorem}\label{FunctorK} For a hypergraph $H \in \mathscr{H}$, define $\mathcal{K}(H):=\{\mathcal{K}(H)_n\}_{n \in \mathbb{Z}_{\geq 0}}$ by \[\mathcal{K}(H)_n := \coprod_{x \in E(H)\coprod V(H)} \Delta[\tilde{f}_H(x)]_n / \sim \] where $\sim$ is the equivalence relation generated by     
\begin{equation}\label{equivalencerelation}
\Delta[\tilde{f}_H(e)]_n \ni (v_{i_0}, \cdots, v_{i_n})_{\tilde{f}_H(x)} \sim (v_{i_0}, \cdots, v_{i_n})_{\tilde{f}_H(x')} \in \Delta[\tilde{f}_H(v)]_n
\end{equation}
for any $x,x' \in E(H)\coprod V(H)$ and $v_{i_0}, \cdots, v_{i_n} \in f_H(x) \cap f_H(x')$. Then $\mathcal{K}(H)$ is a finite symmetric simplicial set. Moreover, $\mathcal{K} : \mathscr{H} \to [!\Delta^{\mathrm{op}}, \mathbf{FinSet}]$ is a functor.
\end{theorem}
\begin{proof}
We denote an equivalence class of $(x, (\underline{v_i})_{\tilde{f}_H(x)})$ in $\mathcal{K}(H)_n$ by $[\underline{v_i}]_x$. For $(\underline{\mu})_{[n]} \in !\Delta([m],[n])$, $\mathcal{K}(H)((\underline{\mu})_{[n]}) : \mathcal{K}(H)_n \to \mathcal{K}(H)_m$ is defined by \[\mathcal{K}(H)((\underline{\mu})_{[n]})([v_{i_0}, \cdots, v_{i_n}]_x) := [v_{i_{\mu_0}}, \cdots, v_{i_{\mu_m}}]_x.\] For any $(\underline{\mu})_{[n]} \in !\Delta([m],[n]), (\underline{\nu})_{[p]} \in !\Delta([n],[p])$ and 
$[\underline{v_i}]_x \in \mathcal{K}(H)_p$, 
\begin{align*}
\mathcal{K}(H)((\underline{\mu})_{[n]}) \circ \mathcal{K}(H)((\underline{\nu})_{[p]})([\underline{v_i}]_x)&=\mathcal{K}(H)((\underline{\mu})_{[n]})([\underline{v_{i_{\nu}}}]_x)\\
&=[\underline{v_{i_{\nu_{\mu}}}}]_x=\mathcal{K}(H)((\underline{\nu})_{[p]} \circ (\underline{\mu})_{[n]}).    
\end{align*}
Hence $\mathcal{K}(H) \in [!\Delta^{\mathrm{op}}, \mathbf{FinSet}]$. \\
Given $\eta= (b,a)\in \mathscr{H}(H, H')$, $\mathcal{K}(\eta)_n : \mathcal{K}(H)_n \to \mathcal{K}(H')_n$ is defined by
\begin{equation}\label{Ksendsmorphism}
\mathcal{K}(\eta)_n\left( [\underline{v_i}]_x \right):=\begin{cases}
    [\underline{a(v_i)}]_{b(e)} & \textrm{if $x=e \in E(H)$} \\
    [\underline{a(v_i)}]_{a(v)} & \textrm{if $x=v \in V(H)$.}
\end{cases}
\end{equation}
For any $(\underline{\mu})_{[n]} \in !\Delta([m],[n])$, $[\underline{v_i}]_x \in \mathcal{K}(H)_n$, 
\begin{align*}
\left(\mathcal{K}(H')((\underline{\mu})_{[n]}) \circ \mathcal{K}(\eta)_n \right)\left([\underline{v_i}]_x \right)&=\left(\mathcal{K}(\eta)_m \circ \mathcal{K}(H)((\underline{\mu})_{[n]}) \right)\left( [\underline{v_i}]_x \right)\\
&=\begin{cases}
    [\underline{a(v_{i_{\mu}})}]_{b(e)} & \textrm{if $x=e \in E(H)$} \\
    [\underline{a(v_{i_{\mu}})}]_{a(v)} & \textrm{if $x=v \in V(H)$.}
\end{cases}
\end{align*}
This implies $\mathcal{K}(\eta) \in [!\Delta^{\mathrm{op}}, \mathbf{FinSet}](\mathcal{K}(H), \mathcal{K}(H'))$. Suppose $\eta \in \mathscr{H}(H_0, H_1)$, $\eta' \in \mathscr{H}(H_1, H_2)$. Then  $\mathcal{K}\left(\eta' \circ \eta \right) = \mathcal{K}\left(\eta' \right) \circ \mathcal{K}\left(\eta \right)$ by Equation \ref{Ksendsmorphism}. Therefore, $\mathcal{K}$ is a functor.
\end{proof}

We say $\mathcal{K}(H)$ the finite symmetric simplicial set induced by hypergraph $H$. There are various operations on $\mathcal{K}(H)$.
\begin{definition}Let $H \in \mathscr{H}$ and $m,n \in \mathbb{Z}_{\geq 0}$.
\begin{itemize}
    \item $p_l : \mathcal{K}(H)_m \to \mathcal{K}(H)_0$ is defined by $p_l := \mathcal{K}(H)((l)_{[m]})$ for $l \in [m]$. Explicitly, it is given by \[p_l\left([v_{i_0}, \cdots, v_{i_m}]_x\right):= [v_{i_l}]_x.\]
    \item The face map $d_l : \mathcal{K}(H)_m \to \mathcal{K}(H)_{m-1}$ is given by    
    \[d_l\left([v_{i_0}, \cdots, v_{i_m}]_x\right)= [v_{i_0}, \cdots, \widehat{v}_{i_l}, \cdots, v_{i_m}]_x.\]
    \item $\vee_l : \mathcal{K}(H)_m \times \mathcal{K}(H)_n \to \mathcal{K}(H)_{m+n-1}$ is defined by \[[\underline{v_i}]_x \vee_l [\underline{v_j}]_x := [v_{i_0}, \cdots, v_{i_{l-1}},v_{j_0}, \cdots, v_{j_n}, v_{i_l}, \cdots, v_{i_m}]_x\] for $l \in [m+1]$.
    \item $\mathfrak{S}_m$ acts on $\mathcal{K}(H)_m$ by $g \cdot := \mathcal{K}(H)(g)$ for $g \in \mathfrak{S}_m$. Explicitly, it is given by \[g \cdot [v_{i_0}, \cdots, v_{i_m}]_x := [v_{i_{g(0)}}, \cdots, v_{i_{g(m)}}]_x.\]
\end{itemize}
\end{definition}
We characterize the preorder on the set of simplices of $\widehat{\mathcal{K}(H)}$ to show that $\mathcal{K}(H)$ is closed and \v{C}ech.

\begin{lemma}\label{preorderonhypergraph}For $H \in \mathscr{H}$ and $[\underline{v_i}]_x \in \mathcal{K}(H)_m, [\underline{v_j}]_{x'} \in \mathcal{K}(H)_n$, the following are equivalent on $\widehat{\mathcal{K}(H)}$.
\begin{enumerate}[label=\textup{(\arabic*)}]
    \item $[\underline{v_i}]_x \lesssim [\underline{v_j}]_{x'}$
    \item $\{\underline{v_i}\} \subseteq \{\underline{v_j} \}.$
\end{enumerate}
\end{lemma}
\begin{proof}
(1)$\implies$(2) Suppose $[\underline{v_i}]_x \lesssim [\underline{v_j}]_{x'}$. There exists $(\underline{\mu})_{[n]} \in !\Delta([m], [n])$ satisfying $K(H)((\underline{\mu})_{[n]})([\underline{v_j}]_{x'}) = [\underline{v_{j_{\mu}}}]_{x'} = [\underline{v_i}]_x$, so $\{\underline{v_i}\} = \{\underline{v_{j_{\mu}}}\} \subseteq \{\underline{v_j}\}$.\\
(2)$\implies$(1) Suppose $\{\underline{v_i}\} \subseteq \{\underline{v_j} \}$. For $k \in [m]$, we can choose $\mu_k \in [n]$ satisfying $v_{j_{\mu_k}}=v_{i_k}$ and define $(\underline{\mu})_{[n]} \in !\Delta([m], [n])$ by $(\underline{\mu})_{[n]}(k) := l_k$. Then $\mathcal{K}(H)((\underline{\mu})_{[n]})([\underline{v_j}]_{x'}) = [\underline{v_{j_{\mu}}}]_{x'} = [\underline{v_i}]_{x'}=[\underline{v_i}]_{x}$ and $[\underline{v_i}]_x \lesssim [\underline{v_j}]_{x'}$.
\end{proof}

\begin{example} The hypergraph $H$ in Example \ref{examplehypergraph}.(2) induces a set of simplices $\widehat{\mathcal{K}(H)}$. See Figure \ref{figurehypergraph}.(b) for some elements of $\widehat{\mathcal{K}(H)}$. 
\begin{itemize}
    \item $[v_2]_{e} \lesssim [v_2, v_3]_e \lesssim [v_0, v_1, v_2]_e$.
    \item $[v_4, v_5]_{e'} \lesssim [v_5, v_4]_{e'}$ and $[v_5, v_4]_{e'} \lesssim [v_4, v_5]_{e'}$. This implies $\lesssim$ is not a partial order since $[v_4, v_5]_{e'} \neq [v_5, v_4]_{e'}$.
    \item $[v_2, v_3]_e=[v_2, v_3]_{e'}$.
\end{itemize}
\end{example}

\begin{proposition}\label{closedunderfiniteintersection}Suppose $H \in \mathscr{H}$ is a hypergraph, $x, x' \in E(H) \coprod V(H)$ and $I \subset \tilde{f}_H(x)$, $I' \subset \tilde{f}_H(x')$.
\begin{enumerate}[label=\textup{(\arabic*)}]
    \item $W_I := U_{[\underline{v_i}]_x} \subset \widehat{\mathcal{K}(H)}$ for any $[\underline{v_i}]_x$ satisfying $\{\underline{v_i}\}=I$ is well-defined.
    \item $W_I \cap W_{I'} = \begin{cases}
    W_{I \cup I'} \textrm{ if there is $x \in E(H)\coprod V(H)$ satisfying $I \cup I' \subset \tilde{f}_H(x)$} \\
    \emptyset \textrm{ otherwise.} 
    \end{cases}$
\end{enumerate}
\end{proposition}
\begin{proof}(1) Suppose $[\underline{v_i}]_x = [\underline{v_j}]_{x'}$. Since $\{\underline{v_i}\} = \{\underline{v_j} \}$ by Equation \ref{equivalencerelation}, $[\underline{v_i}]_x \lesssim [\underline{v_j}]_{x'}$ and $[\underline{v_j}]_{x'} \lesssim [\underline{v_i}]_{x}$. Hence $U_{[\underline{v_i}]_x}=U_{[\underline{v_{j}}]_{x'}}$ and $W_I$ is well-defined.\\
(2) Choose $[\underline{v_i}]_x, [\underline{v_{i'}}]_{x'}$ satisfying $\{\underline{v_i}\}=I$ and $\{\underline{v_{i'}}\}=I'$. Then
\begin{align*}
W_I \cap W_{I'}=U_{[\underline{v_i}]_x} \cap U_{[\underline{v_{i'}}]_{x'}}&=\{[\underline{v_u}]_{x''} \mid  [\underline{v_u}]_{x''} \gtrsim [\underline{v_i}]_x, [\underline{v_u}]_{x''} \gtrsim [\underline{v_{i'}}]_{x'} \} \\
&= \{[\underline{v_u}]_{x''} \mid  \{\underline{v_u}\} \supseteq I, \{\underline{v_u}\} \supseteq I' \}(\because \textrm{Lemma \ref{preorderonhypergraph}}) \\
&=\{[\underline{v_u}]_{x''} \mid  \{\underline{v_u}\} \supseteq I \cup I' \} \\
&=W_{I \cup I'}. 
\end{align*}
If there exists at least one $[\underline{v_u}]_{x''} \in \widehat{\mathcal{K}(H)}$ such that $\{\underline{v_u}\} \supseteq I \cup I'$, $\{[\underline{v_u}]_{x''} \mid  \{\underline{v_u}\} \supseteq I \cup I' \} = W_{I \cup I'}$. Otherwise, $\{[\underline{v_u}]_{x''} \mid  \{\underline{v_u}\} \supseteq I \cup I' \} = \emptyset$.
\end{proof}

Proposition \ref{closedunderfiniteintersection} is used to show that $\mathcal{K}(H)$ is closed, \v{C}ech.

\begin{theorem}For a hypergraph $H \in \mathscr{H}$, $\mathcal{K}(H)$ is closed, \v{C}ech.
\end{theorem}

\begin{proof}
(1) Suppose $[\underline{v_i}]_x, [\underline{v_j}]_{x'} \in \widehat{\mathcal{K}(H)}$. If there is $x'' \in E(H) \coprod V(H)$ satisfying $\{\underline{v_i} \} \cup \{\underline{v_j} \} \subset \tilde{f}_H(x'')$, $U_{[\underline{v_i}]_x} \cap U_{[\underline{v_j}]_{x'}}=W_{\{\underline{v_i} \} \cup \{\underline{v_j} \}}=U_{[\underline{v_i}]_{x''} \vee_0 [\underline{v_j}]_{x''}}$ for $[\underline{v_i}]_{x''} \vee_0 [\underline{v_j}]_{x''} \in \widehat{\mathcal{K}(H)}$ by Proposition \ref{closedunderfiniteintersection}. Otherwise, $U_{[\underline{v_i}]_x} \cap U_{[\underline{v_j}]_{x'}}= \emptyset$ by Proposition \ref{closedunderfiniteintersection}. Hence $\mathcal{K}(H)$ is closed. \\
(2) Proposition \ref{closedunderfiniteintersection} also implies
\begin{align*}
\check{C}\left(\mathcal{K}(H) \right)_n &= \{([v_0]_{v_0}, \cdots, [v_n]_{v_n}) \mid U_{[v_0]_{v_0}}\cap \cdots \cap U_{[v_n]_{v_n}} \neq \emptyset  \} \\
&=\{([v_0]_{x}, \cdots, [v_n]_{x}) \mid v_0, \cdots, v_n \in \tilde{f}_H(x) \textrm{ for some $x \in E(H) \coprod V(H)$} \}
\end{align*}
for any $n \in \mathbb{Z}_{> 0}$.
$\psi_n : \mathcal{K}(H)_n \to \check{C}\left(\mathcal{K}(H) \right)_{n}$ in Definition \ref{Cechdefn} is given by \[\psi_n([v_0, \cdots, v_n]_x) := ([v_0]_x, \cdots, [v_n]_x)\] and it has the inverse $\varphi_n$ with \[\varphi_n\left(([v_0]_{x}, \cdots, [v_n]_{x}) \right) :=[v_0, \cdots, v_n]_x.\]
For any $[v_{i_0}, \cdots, v_{i_n}]_x \in \mathcal{K}(H)_n, ([w_{j_0}]_{x}, \cdots, [w_{j_n}]_{x}) \in \check{C}(\mathcal{K}(H))_n$ and $\mu : [m] \to [n]$,
\begin{align*}
\left(\check{C}(\mathcal{K}(H))(\mu) \circ \psi_n\right)([v_{i_0}, \cdots, v_{i_n}]_x ) &= \check{C}(\mathcal{K}(H))(\mu)\left(([v_0]_x, \cdots, [v_n]_x) \right)\\
&= ([v_{\mu(0)}]_x, \cdots, [v_{\mu(m)}]_x) \\
&= \psi_m([v_{\mu(0)}, \cdots, v_{\mu(m)}]_x) \\
&= \left(\psi_m \circ \mathcal{K}(H)(\mu)\right)([v_{i_0}, \cdots, v_{i_n}]_x)
\end{align*}
and
\begin{align*}
\left(\mathcal{K}(H)(\mu) \circ \varphi_n \right)(([w_{j_0}]_{x}, \cdots, [w_{j_n}]_{x})) &= \mathcal{K}(H)(\mu)([w_{j_0}, \cdots, w_{j_n}]_x)\\
&=[w_{j_{\mu(0)}}, \cdots, w_{j_{\mu(m)}}]_x \\
&= \varphi_m(([w_{j_{\mu(0)}}]_{x}, \cdots, [w_{j_{\mu(m)}}]_{x})) \\
&= \left(\varphi_m \circ \check{C}(\mathcal{K}(H))(\mu) \right)(([w_{j_0}]_{x}, \cdots, [w_{j_n}]_{x})).
\end{align*}
Therefore, $\psi=\{\psi_n\}_{n \in \mathbb{Z}_{\geq 0}} \in [!\Delta^{\mathrm{op}}, \mathbf{FinSet}](\mathcal{K}(H), \check{C}(\mathcal{K}(H)))$ is isomorphism.
\end{proof}

\subsection{Compatibility of cellular sheaf cochain complexes}We show that cellular sheaf cochain complex of an ordered finite abstract simplicial complex $L$ in Equation \ref{CSCCoforderedscpx}, \ref{coboundarymapoforderedscpx} is exactly ordered cellular sheaf cochain complex of $\widehat{\mathcal{K}(L)}$.

\begin{theorem}Suppose $L$ is an ordered finite abstract simplicial complex and $\mathscr{A}$ is an abelian category. For an $\mathscr{A}$-valued cellular sheaf $F$ on $L$, define $F_L \in \mathrm{Cell}(\widehat{\mathcal{K}(L)}, \mathscr{A})$ by $F_L([\underline{v_i}]_x) := F\left( (\underline{v_i}) \right)$. Then $\mathcal{K}(L)_0$ is totally ordered and \[(C^k_F(L, F), \delta^k_F) = (C^k_{\mathrm{ord}, F_L}(\widehat{\mathcal{K}(L)}, F_L), \delta^k_{F_L})\]
for any $k \in \mathbb{Z}_{\geq 0}$.
\end{theorem}

\begin{proof}
A function $\phi : \mathcal{K}(L)_0 \to L_0$ defined by $\phi([v]_v) := v$ is a set isomorphism, so $\mathcal{K}(L)_0$ is also totally ordered preserved by $\phi$. Computations show that
\begin{align*}
C^k_{\mathrm{ord}, F_L}(\widehat{\mathcal{K}(L)}, F_L) &= \bigoplus_{\substack{v_{i_0} < \cdots < v_{i_k} \\ v_{i_0}, \cdots, v_{i_k} \in \tilde{f}_L(x)}}F_L\left(\psi_k^{-1}([v_{i_0}, \cdots, v_{i_k}]_x)\right) \\
&= \bigoplus_{\substack{[v_{i_0}]_x < \cdots < [v_{i_k}]_x \\ v_{i_0}, \cdots, v_{i_k} \in \tilde{f}_L(x)}}F_L\left(([v_{i_0}]_x, \cdots, [v_{i_k}]_x)\right) \\
&= \bigoplus_{\substack{v_{i_0} < \cdots < v_{i_k} \\ v_{i_0}, \cdots, v_{i_k} \in \tilde{f}_L(x)}}F\left((v_{i_0}, \cdots, v_{i_k})\right) \\
&=\bigoplus_{\substack{v_{i_0} < \cdots < v_{i_k} \\ (v_{i_0}, \cdots, v_{i_k}) \in L_k}}F\left((v_{i_0}, \cdots, v_{i_k})\right)=C^k_F(L, F)
\end{align*}
and the restriction of $\delta^k_{F_L}$ to $C^k_{\mathrm{ord}, F_L}(\widehat{\mathcal{K}(L)}, F_L)$ is given by
\begin{align*}
\delta^k_{F_L} &= \underset{\substack{\sigma = [v_{i_0}, \cdots, v_{i_{k+1}}]_x \in \mathcal{K}(L)_{k+1} \\ v_{i_0}< \cdots < v_{i_{k+1}}}}{\bigoplus} \left(\sum_{l \in [k+1]}(-1)^l \cdot (\psi_*F_L)(d_l(\sigma)\lesssim \sigma) \circ \pi_{d_l(\sigma)}\right) \\
&= \underset{\substack{\sigma = (v_{i_0}, \cdots, v_{i_{k+1}}) \in L_{k+1} \\ v_{i_0}< \cdots < v_{i_{k+1}}}}{\bigoplus} \left(\sum_{l \in [k+1]}(-1)^l \cdot F(d_l(\sigma)\lesssim \sigma) \circ \pi_{d_l(\sigma)}\right)\\
&=\delta^k_F.
\end{align*}

\end{proof}

\section{Cellular sheaf Laplacian and its computations}\label{section5} In this section, we provide formulas for degree $k$ cellular sheaf Laplacians on $\widehat{\mathcal{K}(H)}$. 

\begin{proposition}\label{formulaforadjoint} Suppose $H \in \mathscr{H}$ is a hypergraph, $F \in \mathrm{Cell}(\widehat{\mathcal{K}(H)}, \mathbf{Vect}_{\mathbb{R}})$ is a cellular sheaf on $\widehat{\mathcal{K}(H)}$ and $k \in \mathbb{Z}_{\geq 0}$. Given $[\underline{v_i}]_x \in \mathcal{K}(H)_k$, define \[V([\underline{v_i}]_x):= \underset{\{x' \in E(H) \mid \{\underline{v_i} \} \subset f_H(x') \}}{\bigcup}f_H(x')\subset V(H).\]    
\begin{enumerate}[label=\textup{(\arabic*)}]
    \item $\pi_{[\underline{v_i}]_x} \circ (\delta^k_F)^* : C^{k+1}(\widehat{\mathcal{K}(H)},F) \to F([\underline{v_i}]_x)$ is given by \[\underset{\substack{v \in V([\underline{v_i}]_x) \\ l \in [k+1]}}{\sum} (-1)^lF^*([\underline{v_i}]_{x'} \lesssim [\underline{v_i}]_{x'} \vee_l [v]_{x'}) \circ \pi_{[\underline{v_i}]_{x'} \vee_l [v]_{x'}}.\]
    \item $\pi_{[\underline{v_i}]_x} \circ (\delta^k_F)^* : C^{k+1}_{\mathrm{alt}}(\widehat{\mathcal{K}(H)},F) \to F([\underline{v_i}]_x)$ is given by \[(k+2)\left(\sum_{v \in V([\underline{v_i}]_x)} F^*([\underline{v_i}]_{x'} \lesssim [\underline{v_i}]_{x'} \vee_0 [v]_{x'}) \circ \pi_{[\underline{v_i}]_{x'} \vee_0 [v]_{x'}}\right).\]
    \item Suppose $(V(H),<)$ is totally ordered set and $[\underline{v_i}]_x$ satisfies $v_{i_0}< \cdots < v_{i_k}$. For $v \in V([[\underline{v_i}]_x])$, define $l(v) \in [k+1]$ satisfying $v_{i_0} < \cdots < v < v_{i_{l(v)}} < \cdots < v_{i_k}$. Then $\pi_{[\underline{v_i}]_x} \circ (\delta^k_F)^* : C^{k+1}_{\mathrm{ord}}(\widehat{\mathcal{K}(H)},F) \to F([\underline{v_i}]_x)$ is given by \[\sum_{v \in V([\underline{v_i}]_x)} (-1)^{l(v)}F^*([\underline{v_i}]_{x'} \lesssim [\underline{v_i}]_{x'} \vee_{l(v)} [v]_{x'}) \circ \pi_{[\underline{v_i}]_{x'} \vee_{l(v)} [v]_{x'}}.\]
\end{enumerate}
\end{proposition}

\begin{proof}
(1) Equation \ref{cellularcoboundarymap} implies that when $[\underline{v_i}]_x$ is fixed, \[\left(\delta^k_F(s_{[\underline{v_i}]_x})\right)_{[\underline{v_j}]_{x'}} = (-1)^lF([\underline{v_i}]_{x} \lesssim [\underline{v_i}]_{x'} \vee_l [v]_{x'})(s_{[\underline{v_i}]_{x}})\] only if $[\underline{v_j}]_{x'} = [\underline{v_i}]_{x'} \vee_l [v]_{x'}$ for some $v \in V([\underline{v_i}]_x)$, $l \in [k+1]$ and otherwise 0. Hence \[\pi_{[\underline{v_i}]_x} \circ (\delta^k_F)^* = \sum_{\substack{v \in V([\underline{v_i}]_x) \\ l \in [k+1]}} (-1)^lF^*([\underline{v_i}]_{x'} \lesssim [\underline{v_i}]_{x'} \vee_l [v]_{x'}) \circ \pi_{[\underline{v_i}]_{x'} \vee_l [v]_{x'}}.\]
(2) Define $g \in \mathfrak{S}_{k+1}$ by $g(0) := l$, $g(l):=l-1$ for $l \in \{1, \cdots, l\}$ and $g(t):=t$ for $t \in \{l+1, \cdots, k+1\}$. Then $[\underline{v_i}]_{x'} \vee_0 [v]_{x'} = g \cdot \left([\underline{v_i}]_{x'} \vee_l [v]_{x'}\right)$ and \[\pi_{[\underline{v_i}]_{x'} \vee_0 [v]_{x'}} \circ s = \mathrm{sgn}(g) \cdot \pi_{[\underline{v_i}]_{x'} \vee_l [v]_{x'}} \circ s = (-1)^{l} \pi_{[\underline{v_i}]_{x'} \vee_l [v]_{x'}} \circ s\] for $s \in C^{k+1}_{\mathrm{alt}}(\widehat{\mathcal{K}(H)},F)$. Hence
\begin{align*}
\pi_{[\underline{v_i}]_x} \circ (\delta^k_F)^* &= \sum_{\substack{v \in V([\underline{v_i}]_x) \\ l \in [k+1]}} (-1)^lF^*([\underline{v_i}]_{x'} \lesssim [\underline{v_i}]_{x'} \vee_l [v]_{x'}) \circ \pi_{[\underline{v_i}]_{x'} \vee_l [v]_{x'}}\\
&=\sum_{\substack{v \in V([\underline{v_i}]_x) \\ l \in [k+1]}} (-1)^{l+l}F^*([\underline{v_i}]_{x'} \lesssim [\underline{v_i}]_{x'} \vee_0 [v]_{x'}) \circ \pi_{[\underline{v_i}]_{x'} \vee_0 [v]_{x'}}\\
&=(k+2)\left(\sum_{v \in V([\underline{v_i}]_x)} F^*([\underline{v_i}]_{x'} \lesssim [\underline{v_i}]_{x'} \vee_0 [v]_{x'}) \circ \pi_{[\underline{v_i}]_{x'} \vee_0 [v]_{x'}}\right).
\end{align*}
(3) Suppose $[\underline{v_i}]_x \in \mathcal{K}(H)_k, [\underline{v_j}]_{x'} \in \mathcal{K}(H)_{k+1}$ satisfies $v_{i_0}< \cdots < v_{i_k}$ and $v_{j_0} < \cdots < v_{j_{k+1}}$. Then 
\[\left(\delta^k_F(s_{[\underline{v_i}]_x})\right)_{[\underline{v_j}]_{x'}} = (-1)^{l(v)}F([\underline{v_i}]_{x} \lesssim [\underline{v_i}]_{x'} \vee_{l(v)} [v]_{x'})(s_{[\underline{v_i}]_{x}})\] only if $[\underline{v_j}]_{x'} = [\underline{v_i}]_{x'} \vee_{l(v)} [v]_{x'}$ for some $v \in V([\underline{v_i}]_x)$ and otherwise 0. Hence \[\pi_{[\underline{v_i}]_x} \circ (\delta^k_F)^* = \sum_{\substack{v \in V([\underline{v_i}]_x)}} (-1)^{l(v)}F^*([\underline{v_i}]_{x'} \lesssim [\underline{v_i}]_{x'} \vee_{l(v)} [v]_{x'}) \circ \pi_{[\underline{v_i}]_{x'} \vee_{l(v)} [v]_{x'}}.\]
\end{proof}

\begin{theorem}\label{Laplacianformula}Suppose $H \in \mathscr{H}$ is a hypergraph, $F \in \mathrm{Cell}(\widehat{\mathcal{K}(H)}, \mathbf{Vect}_{\mathbb{R}})$ is a cellular sheaf on $\widehat{\mathcal{K}(H)}$, $k \in \mathbb{Z}_{\geq 0}$ and $[\underline{v_i}]_x \in \mathcal{K}(H)_k$.
\begin{enumerate}[label=\textup{(\arabic*)}]
\item For $s \in C^k(\widehat{\mathcal{K}(H)}, F)$, $L^k_{F,+}(s)_{[\underline{v_i}]_x}$ is given by
\begin{align*}
\underset{\substack{l,l' \in [k+1] \\ v \in V([\underline{v_i}]_x)}}{\sum}\Bigl(&(-1)^{l+l'}F^*\left([\underline{v_i}]_{x'} \lesssim [\underline{v_i}]_{x'} \vee_l [v]_{x'}\right)\\
&F\left(d_{l'}([\underline{v_i}]_{x'} \vee_l [v]_{x'}) \lesssim [\underline{v_i}]_{x'} \vee_l [v]_{x'}\right)(s_{d_{l'}([\underline{v_i}]_{x'} \vee_l [v]_{x'})})\Bigr)
\end{align*}
and $L^k_{F,-}(s)_{[\underline{v_i}]_x}$ is given by
\begin{align*}
\underset{\substack{l,l' \in [k] \\ v \in V([\underline{v_i}]_x)}}{\sum}\Bigl(&(-1)^{l+l'}F\left(d_l([\underline{v_i}]_{x'}) \lesssim [\underline{v_i}]_{x'}\right)\\
&F^*\left(d_{l}([\underline{v_i}]_{x'}) \lesssim d_{l}([\underline{v_i}]_{x'}) \vee_{l'} [v]_{x'}\right)(s_{d_{l}([\underline{v_i}]_{x'}) \vee_{l'} [v]_{x'}})\Bigr).
\end{align*}

\item For $s \in C^k_{\mathrm{alt}}(\widehat{\mathcal{K}(H)}, F)$, $L^k_{\mathrm{alt},F,+}(s)_{[\underline{v_i}]_x}$ is given by
\begin{align*}
(k+2)\underset{\substack{l \in [k+1] \\ v \in V([\underline{v_i}]_x)}}{\sum}\Bigl(&(-1)^{l}F^*\left([\underline{v_i}]_{x'} \lesssim [\underline{v_i}]_{x'} \vee_0 [v]_{x'}\right)\\
&F\left(d_{l}([\underline{v_i}]_{x'} \vee_0 [v]_{x'}) \lesssim [\underline{v_i}]_{x'} \vee_0 [v]_{x'}\right)(s_{d_{l}([\underline{v_i}]_{x'} \vee_0 [v]_{x'})})\Bigr)
\end{align*}
and $L^k_{\mathrm{alt}, F,-}(s)_{[\underline{v_i}]_x}$ is given by
\begin{align*}
(k+1)\underset{\substack{l \in [k] \\ v \in V([\underline{v_i}]_x)}}{\sum}\Bigl(&(-1)^{l}F\left(d_l([\underline{v_i}]_{x'}) \lesssim [\underline{v_i}]_{x'}\right)\\
&F^*\left(d_{l}([\underline{v_i}]_{x'}) \lesssim d_{l}([\underline{v_i}]_{x'}) \vee_{0} [v]_{x'}\right)(s_{d_{l}([\underline{v_i}]_{x'}) \vee_{0} [v]_{x'}})\Bigr).
\end{align*}

\item Suppose $(V(H), <)$ is totally ordered set. For $s \in C^k_{\mathrm{ord}}(\widehat{\mathcal{K}(H)}, F)$ and $[\underline{v_i}]_x \in \mathcal{K}(H)_k$ satisfying $v_{i_0} < \cdots < v_{i_k}$, $L^k_{\mathrm{ord},F,+}(s)_{[\underline{v_i}]_x}$ is given by
\begin{align*}
\underset{\substack{l \in [k+1] \\ v \in V([\underline{v_i}]_x)}}{\sum}\Bigl(&(-1)^{l+l(v)}F^*\left([\underline{v_i}]_{x'} \lesssim [\underline{v_i}]_{x'} \vee_{l(v)} [v]_{x'}\right)\\
&F\left(d_{l}([\underline{v_i}]_{x'} \vee_{l(v)} [v]_{x'}) \lesssim [\underline{v_i}]_{x'} \vee_{l(v)} [v]_{x'}\right)(s_{d_{l}([\underline{v_i}]_{x'} \vee_{l(v)} [v]_{x'})})\Bigr)
\end{align*}
and $L^k_{\mathrm{ord}, F,-}(s)_{[\underline{v_i}]_x}$ is given by
\begin{align*}
\underset{\substack{l \in [k] \\ v \in V([\underline{v_i}]_x)}}{\sum}\Bigl(&(-1)^{l+l(v)}F\left(d_l([\underline{v_i}]_{x'}) \lesssim [\underline{v_i}]_{x'}\right)\\
&F^*\left(d_{l}([\underline{v_i}]_{x'}) \lesssim d_{l}([\underline{v_i}]_{x'}) \vee_{l(v)} [v]_{x'}\right)(s_{d_{l}([\underline{v_i}]_{x'}) \vee_{l(v)} [v]_{x'}})\Bigr).
\end{align*} 
\end{enumerate}
\end{theorem}

\begin{proof}
(1) Equation \ref{cellularcoboundarymap} and Proposition \ref{formulaforadjoint}.(1) imply that $L^k_{F,+}(s)_{[\underline{v_i}]_x}$ is equal to
\begin{equation*}
\begin{aligned}
&\sum_{v,l} (-1)^l F^*\left([\underline{v_i}]_{x'} \lesssim [\underline{v_i}]_{x'} \vee_l [v]_{x'}\right) \left(\delta^k_F(s)_{[\underline{v_i}]_{x'} \vee_l [v]_{x'}} \right) \\
= &\sum_{v, l} (-1)^l F^*\left([\underline{v_i}]_{x'} \lesssim [\underline{v_i}]_{x'} \vee_l [v]_{x'}\right)\\
&\quad \hspace{1.5cm}\Bigl(\sum_{l'}(-1)^{l'} F(d_{l'}([\underline{v_i}]_{x'} \vee_l [v]_{x'}) \lesssim [\underline{v_i}]_{x'} \vee_l [v]_{x'})(s_{d_{l'}([\underline{v_i}]_{x'} \vee_l [v]_{x'})})\Bigr) \\
= &\sum_{v, l,l'} (-1)^{l+l'} F^*\left([\underline{v_i}]_{x'} \lesssim [\underline{v_i}]_{x'} \vee_l [v]_{x'}\right)\\
&\quad  \hspace{1.5cm}F\left(d_{l'}([\underline{v_i}]_{x'} \vee_l [v]_{x'}) \lesssim [\underline{v_i}]_{x'} \vee_l [v]_{x'}\right)(s_{d_{l'}([\underline{v_i}]_{x'} \vee_l [v]_{x'})}).
\end{aligned}
\end{equation*}
Same computations proves formula for $L^k_{F,-}(s)_{[\underline{v_i}]_x}$. Proposition \ref{formulaforadjoint}.(2), \ref{formulaforadjoint}.(3) for $(\delta^k_F)^*$ on alternating, ordered cellular sheaf cochains prove (2) and (3).
\end{proof}

\bibliography{references}

\end{document}